\documentclass{article}

\usepackage{amsmath,amsthm, amsfonts}
\usepackage{graphicx}
\usepackage{subcaption}

\usepackage[colorinlistoftodos]{todonotes}
\usepackage{enumitem}
\usepackage{amsfonts}
\usepackage{verbatim}
\usepackage{listings}
\usepackage[margin=1.0in]{geometry}
\usepackage{hyperref}
\usepackage{epstopdf}
\usepackage{multirow}
\usepackage[percent]{overpic}

\usepackage{tikz}
\usetikzlibrary{matrix}
\usepackage{array}

\newcommand\undermat[2]{%
  \makebox[0pt][l]{$\smash{\underbrace{\phantom{%
    \begin{matrix}#2\end{matrix}}}_{\text{$#1$}}}$}#2}

\newtheorem{thm}{Theorem}
\newtheorem{rem}{Remark}

\newtheorem{defn}{Definition}

\newcommand{\sspC}{{\mathcal{C}}}
\newcommand{\sspcoef}{{\mathcal{C}}}

\newcommand{\tfinal}{T_{final}}
\newcommand{\Rtol}{\text{Rtol}}
\newcommand{\Atol}{\text{Atol}}

\newcommand{\pt}{\tilde{p}}

\newcommand{\tauhat}{\tilde{\tau}}
\newcommand{\h}{\Delta t }
\newcommand{\y}{u}
\newcommand{\e}{\mathbf e}
\newcommand{\RK}{\text{RK}}
\newcommand{\wt}{\tilde{w}}
\newcommand{\err}[1]{err_{#1}}

\DeclareMathOperator*{\argmin}{arg\,min}
\newcommand{\bt}{\tilde{b}^T}

\graphicspath{{./figs/}{./figs/work_precision/}}

\title{Embedded error estimation and adaptive step-size control for optimal explicit strong
stability preserving Runge--Kutta methods}
\date{}

\author{Sidafa Conde\footnote{corresponding author}~\footnote{Sandia National Laboratories, Albuquerque, NM 87123, \texttt{sconde@sandia.gov}} \and Imre Fekete\footnote{Department of Applied Analysis and Computational Mathematics, E\"{o}tv\"{o}s Lor\'{a}nd University, MTA-ELTE Numerical Analysis and Large Network Research Groups, Hungary, {\texttt{feipaat@cs.elte.hu}}.
} \and John N. Shadid\footnote{Sandia National Laboratories, Albuquerque, NM 87123, \texttt{jnshadid@sandia.gov}}~\footnote{Department of Mathematics and Statistics, University of New Mexico, Albuquerque, NM 87131}}

\begin{document}
\maketitle

\begin{abstract}
We construct a family of embedded pairs for optimal strong stability preserving explicit Runge--Kutta methods of order $2 \leq p \leq 4$ to be used to obtain numerical solution of spatially discretized hyperbolic PDEs. 
In this construction, the goals include non-defective methods, large region of absolute stability, and optimal error measurement as defined in~\cite{DekkerVerwer84,C00}.
The new family of embedded pairs offer the ability for strong stability preserving (SSP) methods to adapt by varying the step-size based on the local error estimation while maintaining their inherent nonlinear stability properties.
Through several numerical experiments, we assess the overall effectiveness in terms of precision versus work while also taking into consideration accuracy and stability.



\end{abstract}

\section{Introduction}


In what follows, we consider an initial value problem~(IVP) of the form

\begin{equation}
\y'(t) = f(t, \y(t)), \qquad \y(t_0) = \y_0
\label{eqn:ivp}
\end{equation}

\noindent for which the numerical solution is expected to be efficiently obtained by an explicit time-stepping method. {Using the method-of-lines~(MOL) approach, {spatial discretization of the time-dependent partial differential equations~(PDEs) gives rise} to a large system of ordinary differential equations~(ODEs)~\eqref{eqn:ivp}.} The numerical solution of this IVP~\eqref{eqn:ivp} at each time step with an explicit $s$-stage Runge--Kutta~(RK) method $\RK(A,b^T)$ is given by
\begin{equation}
\y_{n+1} = \y_n + \h \sum_{j=1}^{s} b_j f(t_n + c_j\h, Y_j)
\label{eqn:erkCombination}
\end{equation}

\noindent and the internal stages are computed as
\begin{equation}
Y_i = \y_n + \h \sum_{j=1}^{i-1} a_{ij}f(t_n + c_j \h, Y_j), \qquad i = 1, \ldots, s
\end{equation}

\noindent where $\y_n$ is an approximation to the solutions of~\eqref{eqn:ivp} at time $t_n = t_0 + n \h$, $A = (a_{ij})$ and $b^T = (b_j)$ are the coefficients of the method.

In general, most of the time integration for numerical solution of ODEs are computed with a single formula and a fixed step-size. This type of approach can be non-optimal if the solution varies rapidly over small subsets of the integration interval and slowly over larger ones~\cite{Ascher98}.
{Using a uniform step-size, such that $\h_n = \h$ is the same at every time step, the local error would vary at each time step since the error depends on the properties of $f$ and its derivatives.
A small constant step-size may help resolve regions with large variation in solution, at the expense of unnecessary computation in the region of less variability.}

In an attempt to minimize the computational cost and obtain the best possible result, it is necessary to use an adaptive method based on automatic step-size selection.
Numerical computation based on adaptive methods varies the step-size $\h_n$ such that the {local} error can be uniformly distributed at each time step $t_n$. Such an approach is akin to using non uniform Chebyshev nodes in polynomial interpolation~\cite{Trefethen96finitedifference}. 
Practical error estimates are necessary to choose the step-size $\h_n$ sufficiently small to obtain the required precision of the numerical solution and ensure $\h_n$ is large enough to avoid unnecessary computational work~\cite{Ascher98,HairerNorsettWanner93,ShampineReichelt97}.

Several adaptive techniques exist for local truncation error estimations, \emph{e.g.} Taylor series and Richardson extrapolation as well as RK~\cite{HairerNorsettWanner93}. 
{Runge--Kutta formulas can be used to control error and step-size at every step and make step rejection less expensive.}
This development considers adaptive technique based on Runge--Kutta~(RK) methods where two RK methods, {called embedded Runge--Kutta pairs}, one more accurate than the other, are appropriately chosen for which the higher-order method is assumed to approximate~(or be very close to)~to the unknown exact solution.
With this assumption, the difference between the two numerical solutions gives a measurement for the local error. This local truncation error estimate is used to adjust the step-size.
{
An additional advantage of embedded Runge--Kutta methods, compared with Taylor series and Richardson extrapolation, is that the re-evaluation of the stage solutions is no longer necessary. The pairs share the same stage computations, i.e., they have the same $A$. This essentially provides the local error estimation at little to no cost.}



{The general $s$-stage explicit RK pair $\RK(A, b^T,\tilde{b}^T)$ of order $p(p-1)$ allows for adaptive step-size control based on local truncation error estimation. We refer to $\tilde{b}^T$ as RK embedded pair. The extended Butcher tableau of the explicit {formula} embedded pair is represented by~\eqref{eqn:EmbeddedButcherTableau}.}

\begin{equation}
\begin{array}{c|ccccc} 
0\\
c_2    & a_{21} &&&&\\
c_3          & a_{31}  & a_{32} &&&\\
\vdots &  \vdots & \vdots & \ddots&&\\
c_s & a_{s1} & a_{s2} & \ldots & a_{s,s-1} &\\
\hline
       & b_1    & b_2   & b_3 &  \dots & b_s\\
       & \bt_1    & \bt_2 & \bt_3 & \dots & \bt_s\\
\end{array}
\label{eqn:EmbeddedButcherTableau}
\end{equation}

As usual, $c = (c_1, c_2, \ldots, c_s)^T$ is given by $c = A \e$ with $\e = (1, \ldots, 1)^T \in \mathbb{R}^s$. The vectors $b^T$, $\bt$ define the coefficients of the $p$-th and $(p-1)$-th order approximations, respectively. Throughout this paper, we assume that local extrapolation is applied. The integration is advanced using the $p$-th higher order approximation.
The embedded pair produces an estimate of the local truncation error at each Runge--Kutta time step $t_n$ to $t_{n+1} = t_n + {\h}_{n}$ as $e_{n+1} = \| \y_{n+1} - \hat{\y}_{n+1} \|$ where $\hat{\y}_{n+1}$ is the approximate solution obtained by the $(p-1)$-th order $\RK(A,\bt)$ method.


Several embedded formulas have been proposed.
Arguably the best known embedded RK methods are the 3(2) pair of Bogacki and Shampine~\cite{bogackiShampine89} and the 4(5) pair of Dormand and Prince~\cite{dormandPrince198167}, respectively MATLAB's ODE23 and ODE45 implementation~\cite{ShampineReichelt97}. Some other well-known embedded pairs such as Merson, Ceschino, and Zonneveld can be found in \cite{Ascher98,HairerNorsettWanner93,Merson1957}.
{Coupled with robust step control strategy, it has been shown that embedded explicit Runge--Kutta technique is an efficient method for numerical solution of non-stiff initial value problems~\cite{dormandPrince198167,hairer1987solving,Macdonald03,Kennedy2003ARS,ShampineReichelt97}.}




In this work we are interested in embedded pairs for the optimal strong stability preserving~(SSP) explicit Runge--Kutta methods. SSP explicit RK methods are extensively used in numerical computation of hyperbolic conservation 
with total variable diminishing~(TVD) spatial discretization~\cite{G98,G15,Ketcheson11,S88}. 
In many hyperbolic PDE applications, the step-size is controlled by monitoring the CFL number, defined by
\begin{align*}
    \nu & = \frac{c_\text{max} \Delta t}{\Delta x}
\end{align*}
where $c_\text{max}$ is the largest wave speed present.  Since
the Lipschitz constant of the spatial discretization is typically
proportional to $c_\text{max}/\Delta x$, most schemes are stable
up to a particular value of $\nu$.  Inherent in this approach to step-size
control, is the assumption that one can integrate or time-step at or near the largest
(linearly- or nonlinearly-) stable step-size and still achieve an acceptable
level of temporal error. {This is reinforced by experience showing that spatial
error usually dominates temporal error in such problems~\cite{Parsani2013}}.
{However making this assumption and relying on this experience is not enough; an approach that attempts to estimate and control error while adaptively selecting time-steps is important to achieve accurate and inexpensive solutions. Additionally when solving nonlinear PDEs} 
in multiple dimensions, the error estimation provided by embedded RK pairs is essentially
free compared to the expensive evaluation of the right-hand-side.

A second motivation for providing error estimators for SSP methods is
that several optimal SSP methods have good general properties
(useful stability regions, small error coefficients, etc.) and are
frequently used even when SSP theory cannot be applied (i.e.,\
when no forward Euler condition holds), or even for non-hyperbolic PDEs.  
In such situations, how to control the time step-size in practice may be less obvious, while control based on error estimation will be even more useful.

The remainder of this paper is structured as follows: in the next subsection, we briefly review previous work on SSP methods, and present the analytical framework that enables us to construct the new family of embedded pairs. In Section~\ref{sec:methoddesign} we construct the embedded pairs analytically and numerically for optimal explicit SSP Runge--Kutta methods. Section~\ref{sec:numericalresults} contains numerical experiments that compare the newly constructed embedded pairs with existing pairs on several test problems, using different step control strategies, to investigate their performance and robustness. Finally, in Section~\ref{sec:discussion}, we summarize our conclusions and outlook.

\subsection{SSP Runge--Kutta methods}
Strong stability preserving (SSP) time discretization methods were designed to ensure nonlinear stability properties in the numerical solution of spatially discretized hyperbolic PDEs. 
Typically after the spatial discretization we obtain a nonlinear system of ODEs
\begin{equation}\label{eq:semi-discrete}
u_t=F(u),
\end{equation}
where $u$ is a vector of approximations to the exact solution of the PDE. We assume that the semi-discretization \eqref{eq:semi-discrete} and a convex functional $||\cdot||$ (or norm, semi-norm) are given, and that there exists a $\Delta t_{\text{FE}}$ such that the forward Euler condition
\begin{align}\label{eq:FEcond}
||u+\Delta tF(u)||\leq ||u||\ \text{for}\ 0\leq \Delta t\leq \Delta t_{\text{FE}}
\end{align}
holds for all $u$.  An explicit Runge--Kutta (ERK) method is called SSP if the estimate
$$
||u_{n+1}||\leq||u_n||
$$
holds for the numerical solution of \eqref{eq:semi-discrete}, whenever \eqref{eq:FEcond} holds and $\Delta t\leq \sspC\Delta t_{\text{FE}}$. The constant $\sspC$ is called the SSP coefficient.

Instead of using ERK methods in the Butcher form, Shu and Osher suggested another representation of ERK methods in order to design high-order time-discretization methods for problem \eqref{eq:semi-discrete} in \cite{SO88}. By using their idea we can see that certain ERK methods can be rewritten as a convex combination of forward Euler methods. This decomposition is a sufficient and necessary condition for the SSP property. For the relation between Shu-Osher and Butcher representations and for a complete introduction to the topic we recommend monograph \cite{Ketcheson11}. A brief summary of current theoretical results can be found in the recent review \cite{G15}. In what follows we use the corresponding Butcher form of explicit SSP Runge--Kutta methods. Next, we highlight results of \cite{G15} and  \cite{Ketcheson11} which will be used in this paper.
\begin{thm}[\cite{Ketcheson11}, Theorem 3.2.]\label{thm:intro2} Let us consider the matrix
\[
K= 
  \begin{pmatrix} A & 0\\
  b^T & 0
 \end{pmatrix}
\]

and the SSP conditions
\begin{subequations} \label{eq:SSPcon}
\begin{align}
	K(I+rK)^{-1}&\geq 0 \label{eq:SSPcondA}\\
rK(I+rK)^{-1}e&\leq e.\label{eq:SSPcondb}
\end{align}
\end{subequations}

Then, the SSP coefficient of the ERK method is 
\begin{align*}
	\sspC(A,b^T) = \sup\left\{r : (I + rK)^{-1} \text{exists and conditions \eqref{eq:SSPcondA}-\eqref{eq:SSPcondb} hold} \right\}.
\end{align*}
\end{thm}
\begin{thm}[\cite{Ketcheson11}, Observation 5.2.]\label{thm:intro3} Consider an ERK method. If the method has positive SSP coefficient $\sspC(A,b^T)$, then $A\geq 0$ and $b^T\geq 0$. 
\end{thm}
There is no ERK method of order $p\geq 5$ with positive SSP coefficient~\cite{G15}. Therefore, we only give the order conditions up to fourth order. These are
\begin{subequations} \label{eq:ordercond}
\begin{align}
    b^Te & = 1, & (p=1) \label{eq:ordercond1} \\
    b^T c & = \frac{1}{2}, & (p=2) \label{eq:ordercond2}  \\
    b^T c^2 & = \frac{1}{3}, & (p=3) \label{eq:ordercond3a} \\
    b^T \left(\frac{c^2}{2!} - Ac\right) & = 0, & (p=3)\label{eq:ordercond3b}\\
    b^Tc^3&=\frac{1}{4}, & (p=4)\label{eq:ordercond4a}\\
    b^TA \left(\frac{c^2}{2!} - Ac\right) & = 0, & (p=4)\label{eq:ordercond4b}\\
    b^T \left(\frac{c^3}{3!} - \frac{Ac^2}{2!}\right) & = 0, & (p=4)\label{eq:ordercond4c}\\
    b^T\text{diag}(c) \left(\frac{c^2}{2!} - Ac\right) & = 0, & (p=4)\label{eq:ordercond4d}
\end{align}
\end{subequations}
where $\text{diag}(c)$ is the square diagonal matrix with the elements of vector $c$ on the main diagonal.

\section{Embedded pairs for optimal explicit SSP Runge--Kutta methods}\label{sec:methoddesign}


We introduce the notation SSPERK$(s,p)$ for optimal explicit SSP Runge--Kutta methods, where $s$ and $p$ 
refer to the number of stages and order, respectively. We give below the desired properties for embedded methods.

\begin{itemize}
\item[(i),] The embedded method is order of $p-1$, i.e., it has one order less than the SSPERK method.
\item[(ii),] The embedded method is non-defective, i.e., it violates all of the $p$-th order conditions. For practical importance see Remark~\ref{remark:nonDefective}.
\item[(iii),] Whenever possible, the embedded method has rational coefficients and a simple structure.
\item[(iv),]  The embedded method has maximum SSP coefficient $\tilde{\sspC}$, where $\tilde{\sspC}$ is the SSP coefficient of the the optimal SSPERK method; if this is not the case, then we are lookign for embedded SSPERK methods with smaller SSP coefficient or simply embedded ERK methods.
\end{itemize}

Taking into account the desired properties (i)-(iv), we seek an embedded pair $\bt$, with the stage coefficient $A$ from a SSPERK method, such that these satisfy the following simplified optimization problem

\begin{align}
 \begin{pmatrix} A & 0\\
  \tilde{b}^T & 0
 \end{pmatrix}\left(I+\tilde{\sspC}\begin{pmatrix} A & 0\\
  \tilde{b}^T & 0
 \end{pmatrix}\right)^{-1}\geq 0, \label{eq:SSPcondAcomponent}\\
  \left|\left|\tilde{\sspC}\begin{pmatrix} A & 0\\
  \tilde{b}^T & 0
 \end{pmatrix}\left(I+\tilde{\sspC}\begin{pmatrix} A & 0\\
  \tilde{b}^T & 0
 \end{pmatrix}\right)^{-1}\right|\right|_{\infty}\leq 1, \label{eq:SSPcondbcomponent}
 \end{align}
 \begin{align} 
\text{the appropriate order conditions}\ \eqref{eq:ordercond1}-\eqref{eq:ordercond4d}\ \text{and property (ii) are fulfilled},&\label{eq:ordercondLP}
\end{align}

\noindent where \eqref{eq:SSPcondAcomponent}-\eqref{eq:SSPcondbcomponent} are equivalent with \eqref{eq:SSPcondA}-\eqref{eq:SSPcondb} and $||\cdot||_{\infty}$ denotes the induced matrix norm. The reason why we called \eqref{eq:SSPcondAcomponent}-\eqref{eq:SSPcondbcomponent} a simplified optimization problem is that we fix the SSP coefficient $\tilde{\sspC}$. Due to Theorem \ref{thm:intro3} and \eqref{eq:ordercond1} we have $0\leq \bt \leq \mathbf{e}$.

The newly constructed pairs satisfy the desired properties (i)-(iv) and have optimum region of absolute stability. 
The absolute stability region is given by
\begin{equation}
S:=\{z\in\mathbb{C}: |\psi(z)|\leq 1\},
\label{eq:regionAS}
\end{equation}
where $\psi(z)$ is the absolute stability function of the given SSPERK method. The following definitions can be found in
~\cite{DekkerVerwer84, K91, M90}.
\begin{defn}\label{defn:LSRAI} The absolute stability real axis inclusion is the radius of the largest
interval on the real axis that is contained in the absolute stability region. Specifically,
$$
\delta_{R}:=\sup\{\gamma:\gamma\geq\ 0\ and\ l(-\gamma,\gamma)\subset S\},
$$
where $l(z_1,z_2)$ is the line segment connecting $z_1,z_2\in\mathbb{C}$.
\end{defn}
\begin{defn}\label{defn:LSIAI} The absolute stability imaginary axis inclusion is the radius of the largest
interval on the imaginary axis that is contained in the absolute stability region. Specifically,
$$
\delta_{I}:=\sup\{\gamma:\gamma\geq\ 0\ and\ l(-i\gamma,i\gamma)\subset S\},
$$
where $l(z_1,z_2)$ is the line segment connecting $z_1,z_2\in\mathbb{C}$.
\end{defn}
\begin{defn}\label{defn:CCR} An ERK method is called circle contractive if for $r>0$ 
\begin{equation}\label{eq:gendisk}
|\psi(z)|\leq 1\ \text{for all}\ z\in D(r).
\end{equation}
The radius of circle contractivity is the radius of the largest generalized disk $D(r)$ for which \eqref{eq:gendisk} holds, where
\begin{alignat}{4}
D(r)=&\{\lambda\in\mathbb{C}: |\lambda+r|\leq r\}.\nonumber 
\end{alignat}
The radius of circle contractivity will be denoted by $\delta_{C}$.
\end{defn}
\begin{defn}\label{defn:AMR} The radius of absolute monotonicity of the stability function $\psi(z)$ is the largest value of $r$ such that $\psi(z)$ and all of its derivatives exist and are non-negative for $z\in (-r,0]$. It will be denoted by $R(\psi)$.
\end{defn}

Throughout the paper, we will simply refer to values of Definitions \ref{defn:LSRAI}-\ref{defn:AMR} as stability radius measurements. Apart from satisfying the desired properties (i)-(iv) listed above, more constraints are imposed on the pairs. We quantify the $L_2$ and $L_{\infty}$ principal errors of the pairs following~\cite{C00}:
\begin{alignat*}{2}
A^{(p+1)}_2&=||\tau^{(p+1)}||_2 &\quad \text{and}\quad \tilde{A}^{(p)}_2&=||\tilde{\tau}^{(p)}||_2, \\
A^{(p+1)}_{\infty}&=||\tau^{(p+1)}||_{\infty}&\quad \text{and}\quad \tilde{A}^{(p)}_{\infty}&=||\tilde{\tau}^{(p)}||_{\infty},
\end{alignat*}
where $\tau^{(p+1)}$ is the error coefficient vector of SSPERK method of order $p$. The vector $\tilde{\tau}^{(p)}$ corresponds to the error coefficient vector for the embedded pair of order $\pt = p -1$. Additional error controls are defined as
\begin{alignat*}{2}
B^{(p+1)}_2&=\displaystyle\frac{A^{(p+1)}_2}{A^{(p)}_2}& \tilde{B}^{(p)}_2&=\displaystyle\frac{\tilde{A}^{(p)}_2}{\tilde{A}^{(p-1)}_2},\\
B^{(p+1)}_{\infty}&=\displaystyle\frac{A^{(p+1)}_{\infty}}{A^{(p)}_{\infty}}& \tilde{B}^{(p)}_{\infty}&=\displaystyle\frac{\tilde{A}^{(p)}_{\infty}}{\tilde{A}^{(p-1)}_{\infty}},\\
    C^{(p+1)}_2&=\displaystyle\frac{\| \tauhat^{(p+1)} - \tau^{(p+1)} \|_2}{\|\tauhat^{(p)}\|_2} & C^{(p+1)}_\infty&=\displaystyle\frac{\| \tauhat^{(p+1)} - \tau^{(p+1)} \|_\infty}{\|\tauhat^{(p)}\|_\infty} \\
    {D}&=\max\{|a_{ij}|,|c_i|,|b_i^T|,|\tilde{b}_i^T|\}.\\
\end{alignat*}

\noindent and impose that the values $\tilde{A}^{(p)}_2$  and $\tilde{A}^{(p+1)}_{\infty}$ should be as small as possible, while the values  $B^{(p+1)}_2, B^{(p+1)}_{\infty}$, $ C^{(p+1)}_2$ and $C^{(p+1)}_{\infty}$ should be close to one. Furthermore, the magnitude of $D$ should be small, too. This is equivalent to Dormand and Prince's constraints on $|c_i| \leq 1$ and avoiding large values of $b_i^T$ and $a_{ij}$ to circumvent considerable rounding errors in practical applications~\cite{dormandPrince198167}. 
Throughout the paper we will simply refer to these values as error measurements. 

\begin{rem}[optimal SSP]\label{remark:optimalSSP}
    Since the method $\RK(A,b^T)$ is the optimal SSP method of order $p$, the quantities $A^{(p+1)}_2$ and $A^{(p+1)}_{\infty}$ are fixed and beyond our control.
\end{rem}

\begin{rem}[non-defective]\label{remark:nonDefective}
    For a defective embedded pair, $\tauhat_i^{p} = 0$, since it satisfies all of the order $p = \pt + 1$ algebraic order conditions. Consequently $C^{(p+1)}_2, C^{(p+1)}_\infty \rightarrow \infty$. Since the higher-order method is the optimal SSP method, we know it's a non-defective method of order $p$. Otherwise, it would be a non-defective method of order $p+1$.
\end{rem}

Many of the optimal SSP methods, such as {SSPERK(5,4)}~\cite{K91,RS02}, optimal implicit SSP RK~\cite{Ketcheson11}, and the newly developed optimized SSP IMEX methods in~\cite{Conde2017}, were obtained numerically following the methodology in~\cite{Ketcheson08}. Similarly, in constructing an efficient and robust embedded pair for these numerically optimal methods, we construct an analogous optimization problem:

\begin{align}
    F(A, b^T, \wt) = \begin{bmatrix} \tilde{A}_2^{(p)} &
    \tilde{A}_\infty^{(p)} &
\left(B_2^{(p+1)} - 1 \right) & 
\left(B_\infty^{(p+1)} - 1\right) &
 \left(C_2^{(p+1)} - 1 \right)&
 \left(C_\infty^{(p+1) }- 1 \right)
   \end{bmatrix}
   \label{eq:numOptimization}
\end{align}

\begin{equation*}
\begin{aligned}
& \argmin_{\wt}
& & \| F(A, b^T, \wt) \|_\infty \\
& \text{subject to}
& & \tau_k(A, \wt) = 0, \; k = 1, \ldots, \pt \\
& & & \wt \geq 0
\end{aligned}
\end{equation*}

\noindent where $\tau_k(A, \wt)$ is the necessary algebraic order conditions for the embedded method of order $\pt$ \eqref{eq:ordercond}. Here, we have used $\wt$ to represent the embedded weight obtained numerically. Using MATLAB's \emph{fmincon} built-in function, we search for $\wt$ that minimizes the cost function best with respect to stability. The positivity constraint on the embedded weights, $0\leq \wt \leq \mathbf{e}$, is  due to Theorem \ref{thm:intro3} and \eqref{eq:ordercond1}. 
Moreover, we find the above constraints always enforce $\| D \| \leq 1$ .
For the optimization problem, we do not impose the method to have rational coefficients since many of the optimal SSP methods do not satisfy this restriction - a requirement not enforced by Carpenter, Kennedy, Lewis in~\cite{Carpenter2000}.
%


\subsection{Embedded pairs for SSPERK(s,2) methods}\label{sec:SSPERK(s,2)}
Gottlieb and Shu made the first step to characterize SSPERK$(s,2)$ methods after they gave the SSPERK$(2,2)$ method~\cite{G98}. Later, Ruuth and Spiteri proved that SSPERK$(s,2)$ methods have $\sspC=s-1$ and gave their Shu-Osher representation~\cite{RS02}. The Butcher form of SSPERK$(s,2)$ methods is given in Table \ref{table:SSPERK(s,2)}.

\begin{table}[ht!]
\begin{center}
\begin{tabular}{c|ccccc} 
0\\
$\frac{1}{s-1}$    & $\frac{1}{s-1}$ &&&&\\
$\frac{2}{s-1}$          & $\frac{1}{s-1}$  & $\frac{1}{s-1}$ &&&\\
$\vdots$ &  $\vdots$ & $\vdots$ & $\ddots$&&\\
$1$ & $\frac{1}{s-1}$ & $\frac{1}{s-1}$ & $\ldots$ & $\frac{1}{s-1}$ &\\
\hline\rule{0pt}{1.0\normalbaselineskip}
& $\frac{1}{s}$ & $\frac{1}{s}$ & $\dots$ & $\frac{1}{s}$ & $\frac{1}{s}$
\end{tabular} \caption{Butcher form of SSPERK$(s,2)$ methods.}\label{table:SSPERK(s,2)}
\end{center}
\end{table}

Taking into account the desired properties (i)-(iv) we are looking for non-defective embedded pairs for SSPERK$(s,1)$ methods with $\tilde{\sspC}=s-1$. Below we present only two  best embedded pairs of our numerical search which also have simple structures. These are the pairs
$$
\tilde{b}_1^T=\left[\frac{1}{s-1},\ldots,\frac{1}{s-1},0\right]
$$
and
$$
\tilde{b}_2^T=\left[\frac{s+1}{s^2},\frac{1}{s},\ldots,\frac{1}{s},\frac{s-1}{s^2}\right].
$$

In order to visually demonstrate the relative stability of the pairs $\bt_1$ and $\bt_2$ related to desired properties (i)-(iv) we plot their stability regions and the stability region of SSPERK$(s,2)$ methods for $s=2,4,6, 8$. Based on just Figure  \ref{fig:StabRegSSPERK(s,2)}, these methods can be considered equally good for $s\geq 4$. However, if we consider the introduced stability radius measurement values and error measurement values from Section \ref{sec:methoddesign}, it turns out that it is worth choosing the embedded pair $\bt_2$. The corresponding stability radius measurement values and error measurement values can be found in~\cite{CondeFekete18}.

\begin{figure}[ht!]
    \centering
    \begin{subfigure}[b]{0.25\textwidth}
        \includegraphics[width=\textwidth,height=\textwidth]{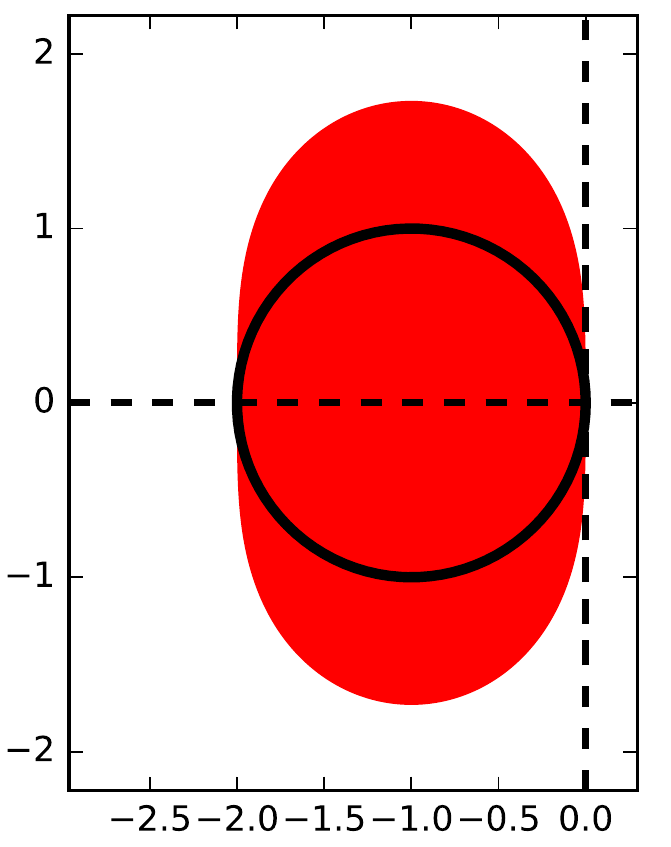} \caption{SSPERK$(2,2)-\bt_1$}
        \label{fig:SSPERK22}
    \end{subfigure}~ 
    \begin{subfigure}[b]{0.25\textwidth}
        \includegraphics[width=\textwidth,height=\textwidth]{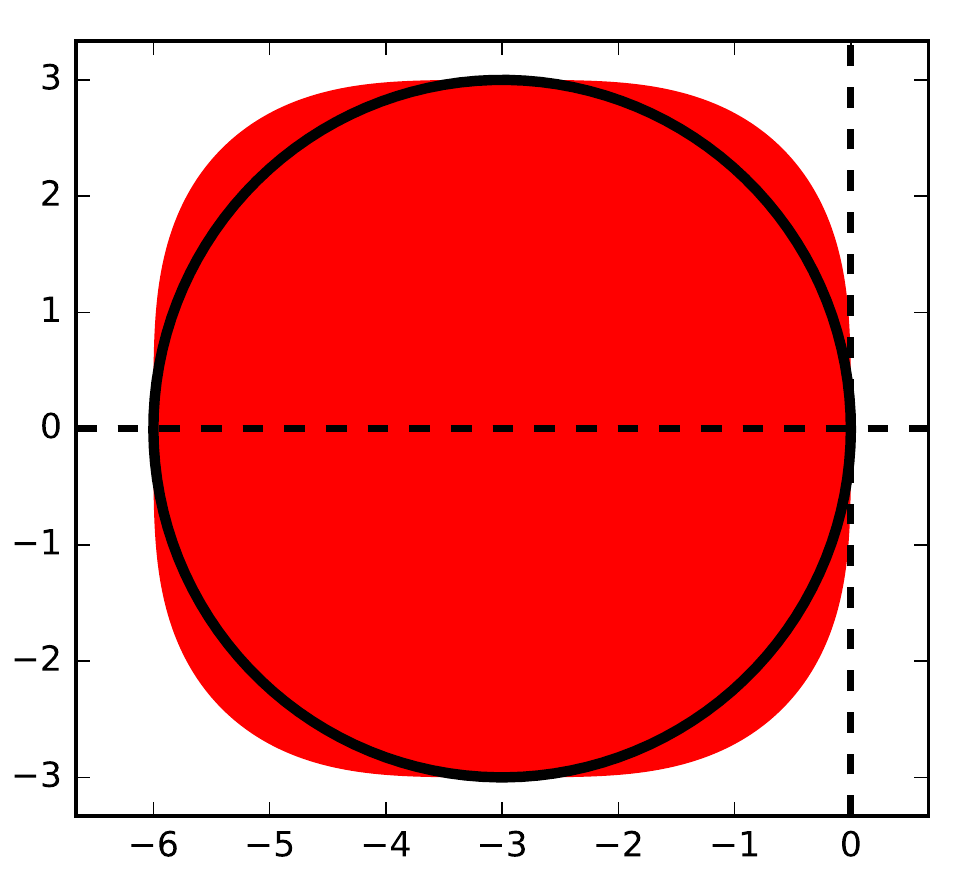}
        \caption{SSPERK$(4,2)-\bt_1$}
        \label{fig:SSPERK42}
    \end{subfigure}~
    \begin{subfigure}[b]{0.25\textwidth}
        \includegraphics[width=\textwidth,height=\textwidth]{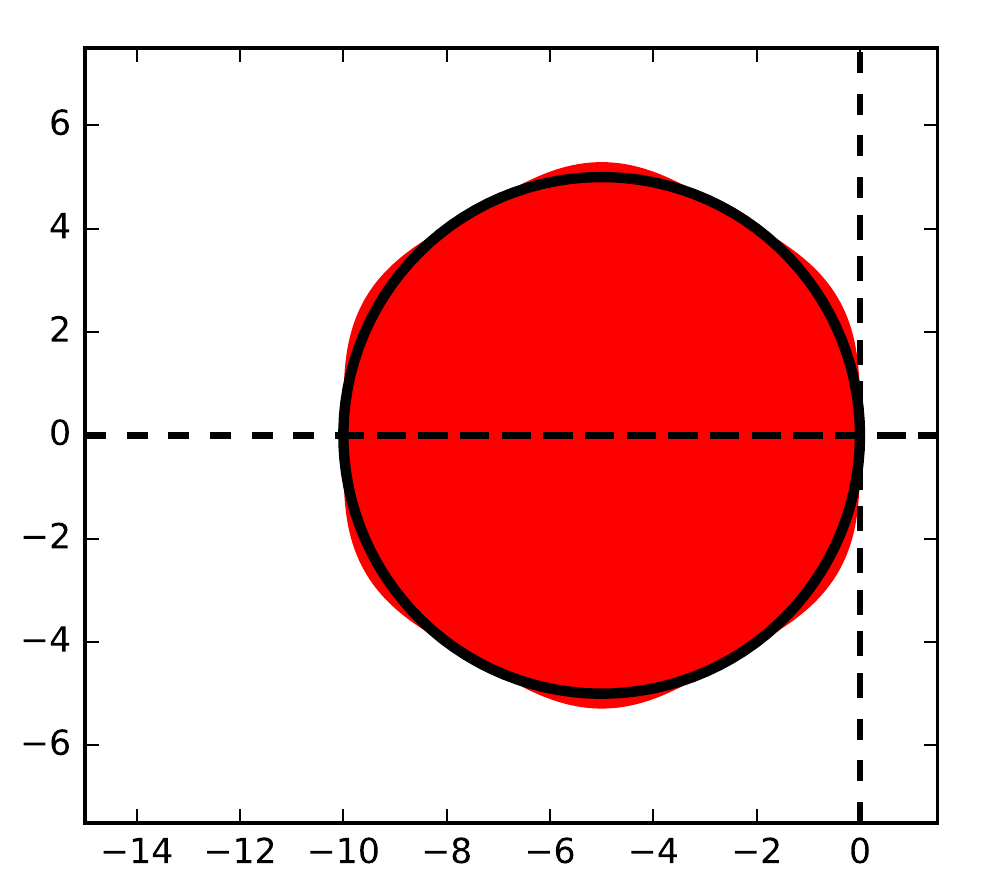}
        \caption{SSPERK$(6,2)-\bt_1$}
        \label{fig:SSPERK62}
    \end{subfigure}~
    \begin{subfigure}[b]{0.25\textwidth}
        \includegraphics[width=\textwidth,height=\textwidth]{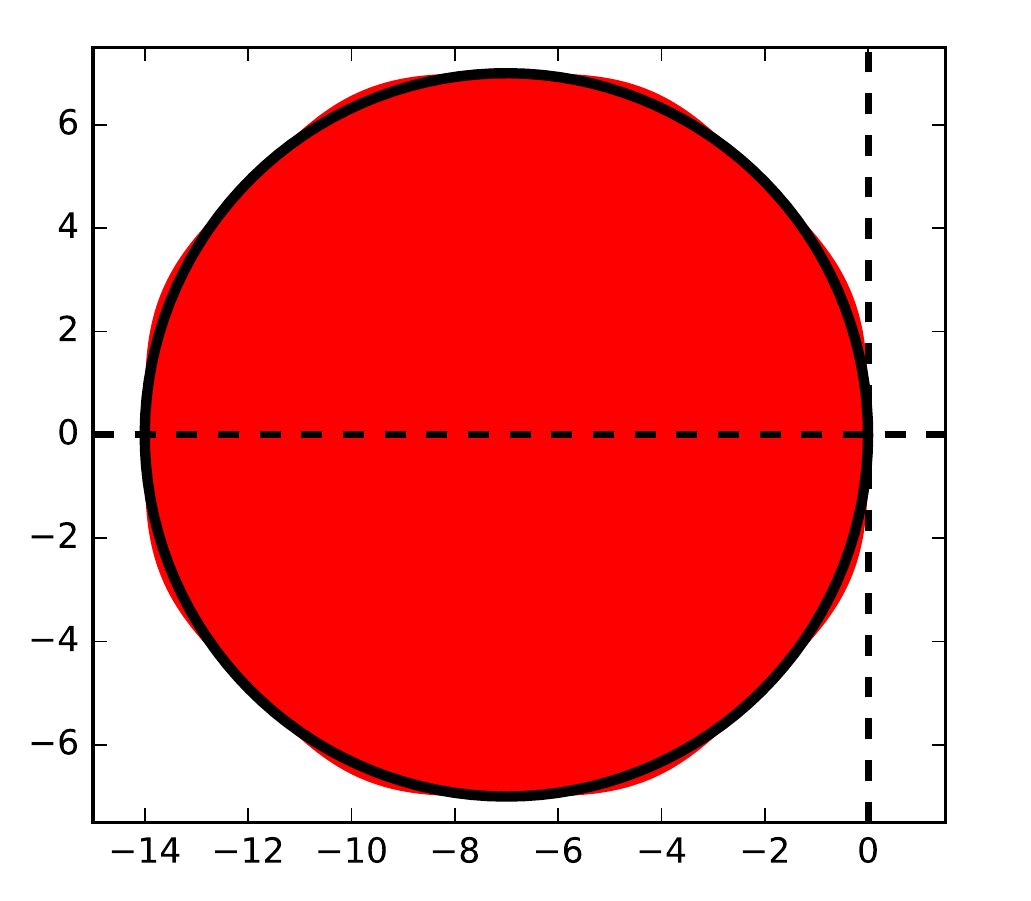}
        \caption{SSPERK$(8,2)-\bt_1$}
        \label{fig:SSPERK82}
    \end{subfigure}\\
        \begin{subfigure}[b]{0.25\textwidth}
        \includegraphics[width=\textwidth,height=\textwidth]{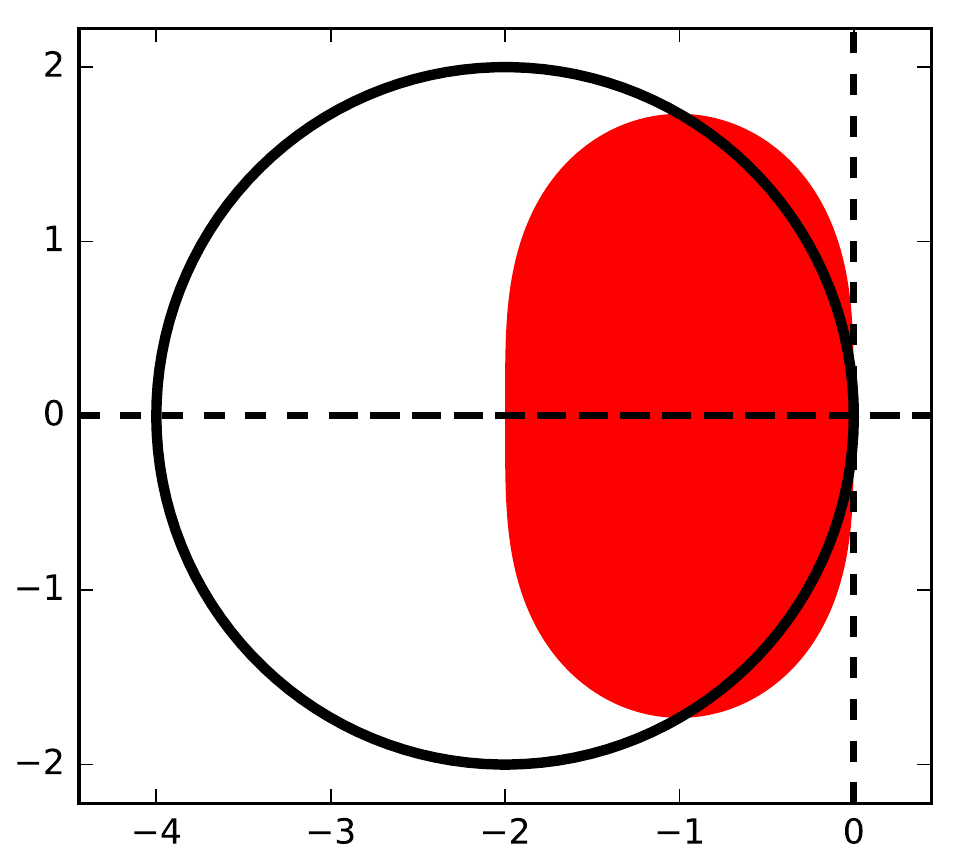} \caption{SSPERK$(2,2)-\bt_2$}
        \label{fig:SSPERK22b}
    \end{subfigure}~ 
    \begin{subfigure}[b]{0.25\textwidth}
        \includegraphics[width=\textwidth,height=\textwidth]{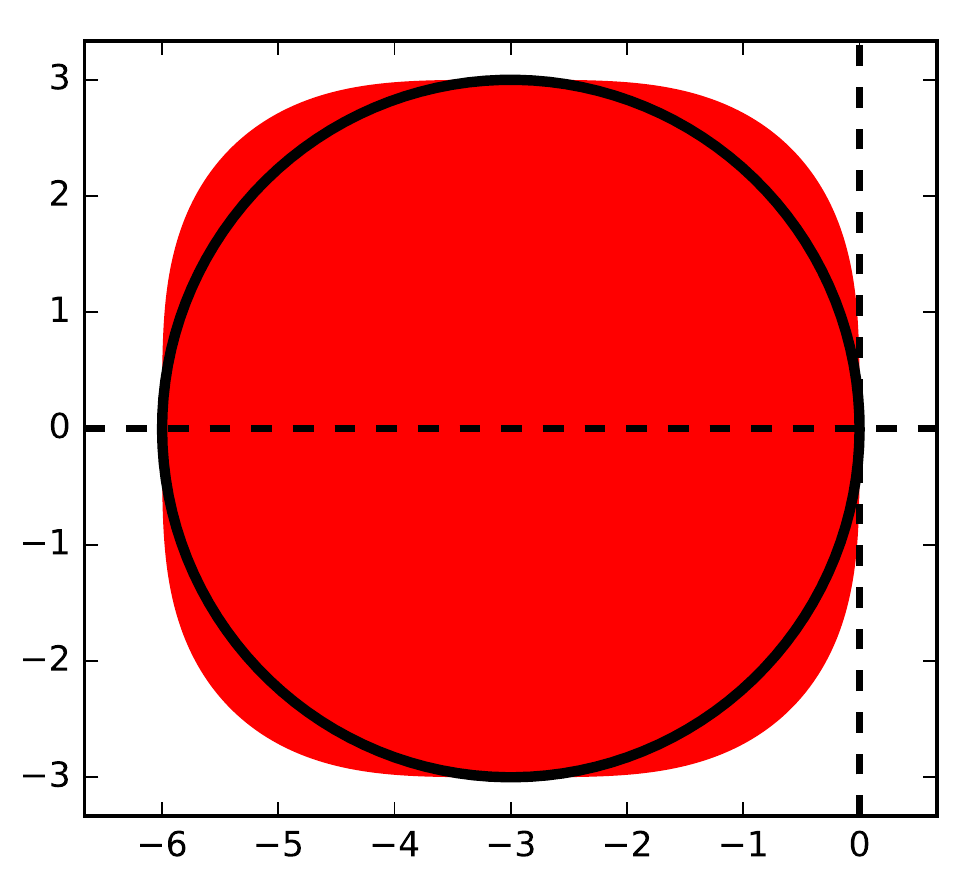}
        \caption{SSPERK$(4,2)-\bt_2$}
        \label{fig:SSPERK42b}
    \end{subfigure}~
    \begin{subfigure}[b]{0.25\textwidth}
        \includegraphics[width=\textwidth,height=\textwidth]{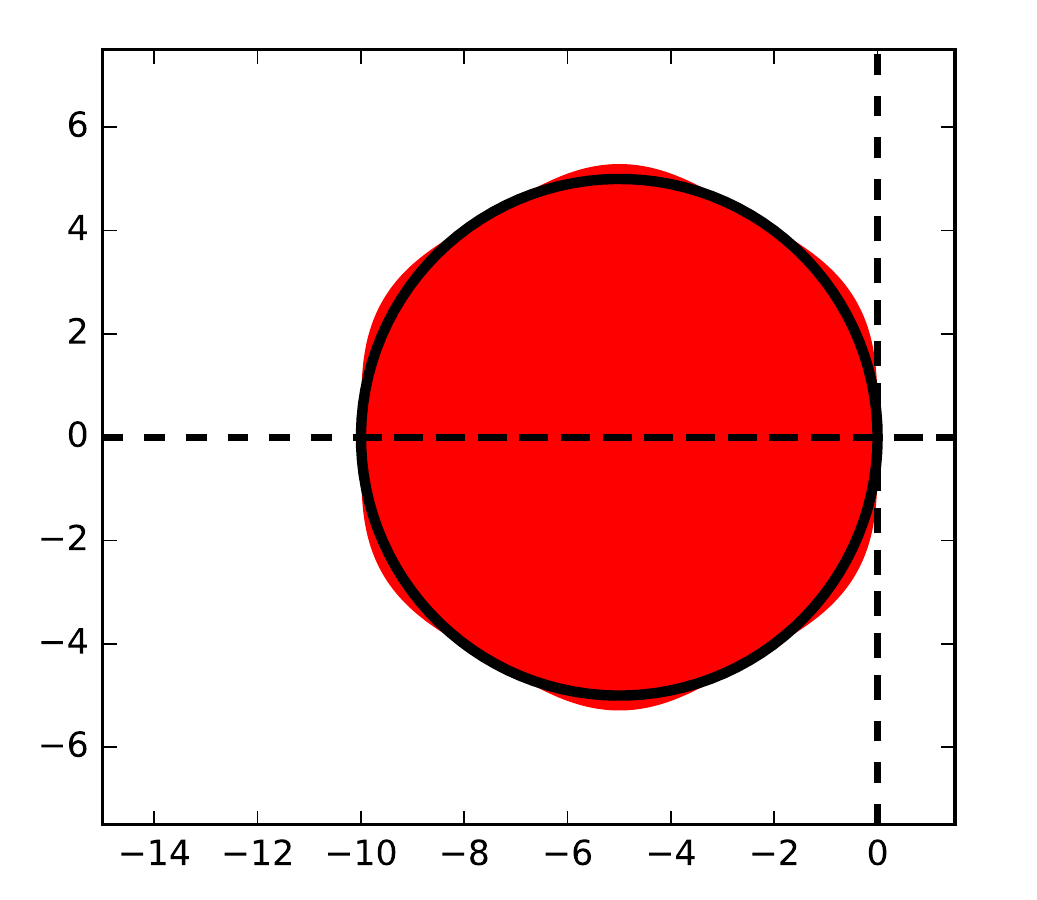}
        \caption{SSPERK$(6,2)-\bt_2$}
        \label{fig:SSPERK62b}
    \end{subfigure}~
    \begin{subfigure}[b]{0.25\textwidth}
        \includegraphics[width=\textwidth,height=\textwidth]{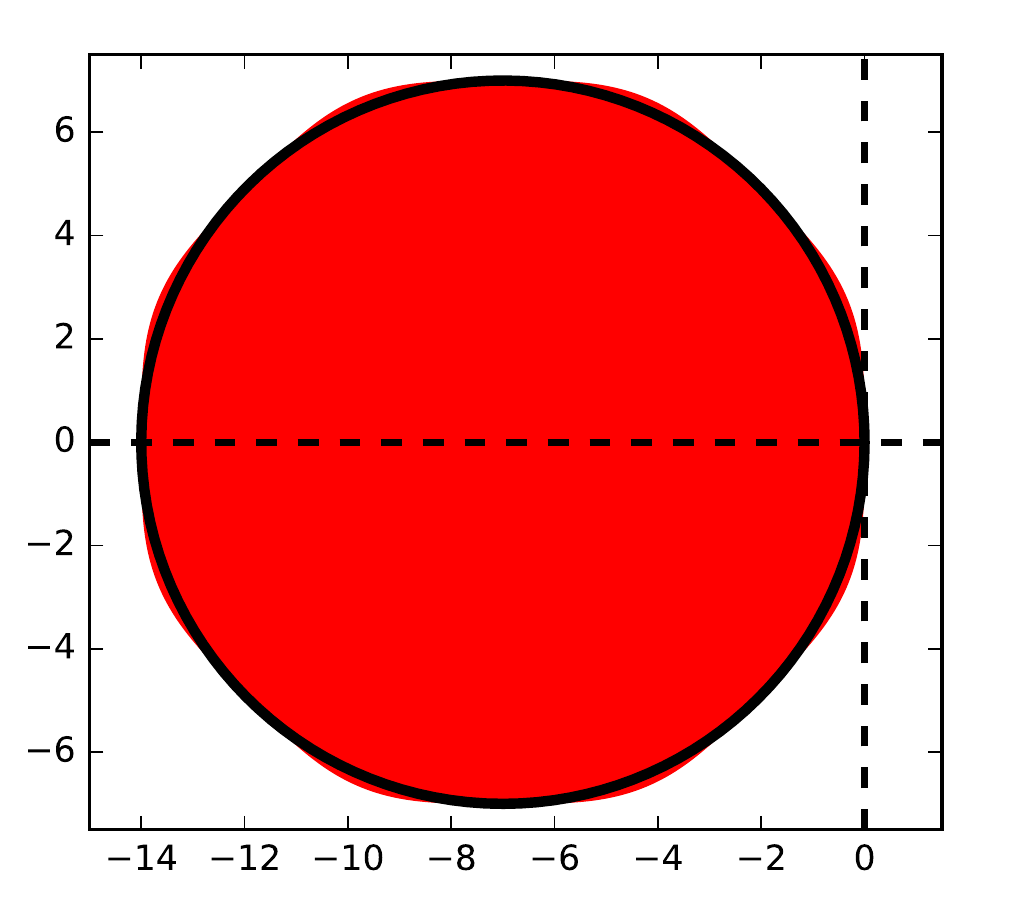}
        \caption{SSPERK$(8,2)-\bt_2$}
        \label{fig:SSPERK82b}
    \end{subfigure}
    \caption{The stability regions of SSPERK$(s,2)$ methods (red) and the black contours of the embedded SSPERK$(s,1)$ methods with $\bt_1$ (top row) and $\bt_2$ (bottom row) for $s=2,4,6,8$. It is more evident when comparing SSPERK$(2,2)-\bt_1$~(Fig. \ref{fig:SSPERK22}) and SSPERK$(2,2)-\bt_2$~(Fig. \ref{fig:SSPERK22b}) that the family of pairs $\bt_2$ are superior as they offer larger stability region and {better error measurements}.}\label{fig:StabRegSSPERK(s,2)}
\end{figure}

Furthermore, two second order pairs obtained from the numerical optimization~(Eq.\eqref{eq:numOptimization}) are SSPERK$(2,2)-\wt$ and SSPERK$(3,2)-\wt$, stability regions in~(Figure-\ref{fig:NumericalStabRegSSPERK(s,2)}).

\begin{figure}[htbp]
    \centering
    \begin{subfigure}[b]{0.47\textwidth}
        \includegraphics[width=\textwidth,height=\textwidth]{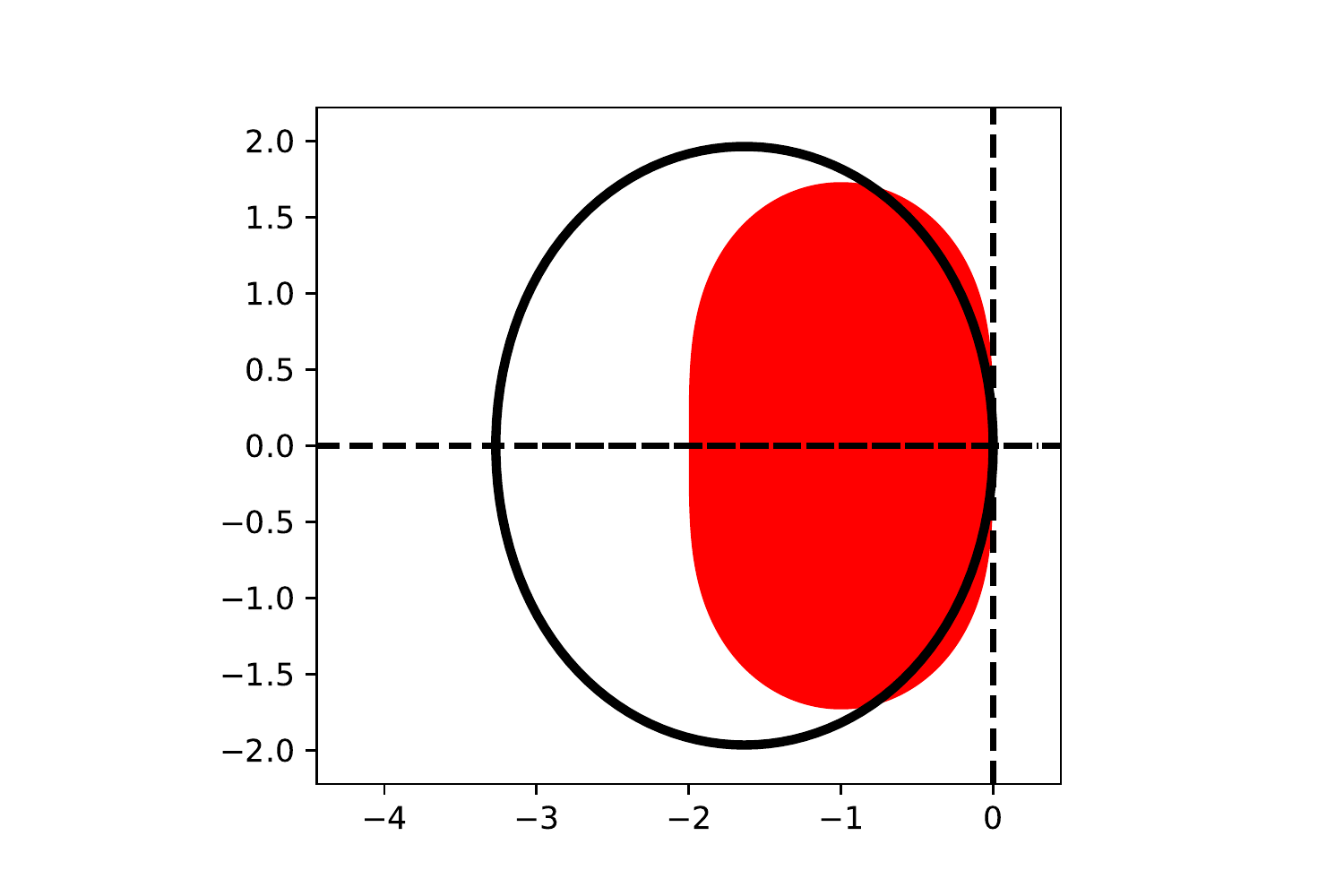} \caption{SSPERK$(2,2)-\wt$}
        \label{fig:SSPERK22BestNew01}
    \end{subfigure}
    \begin{subfigure}[b]{0.47\textwidth}
        \includegraphics[width=\textwidth,height=\textwidth]{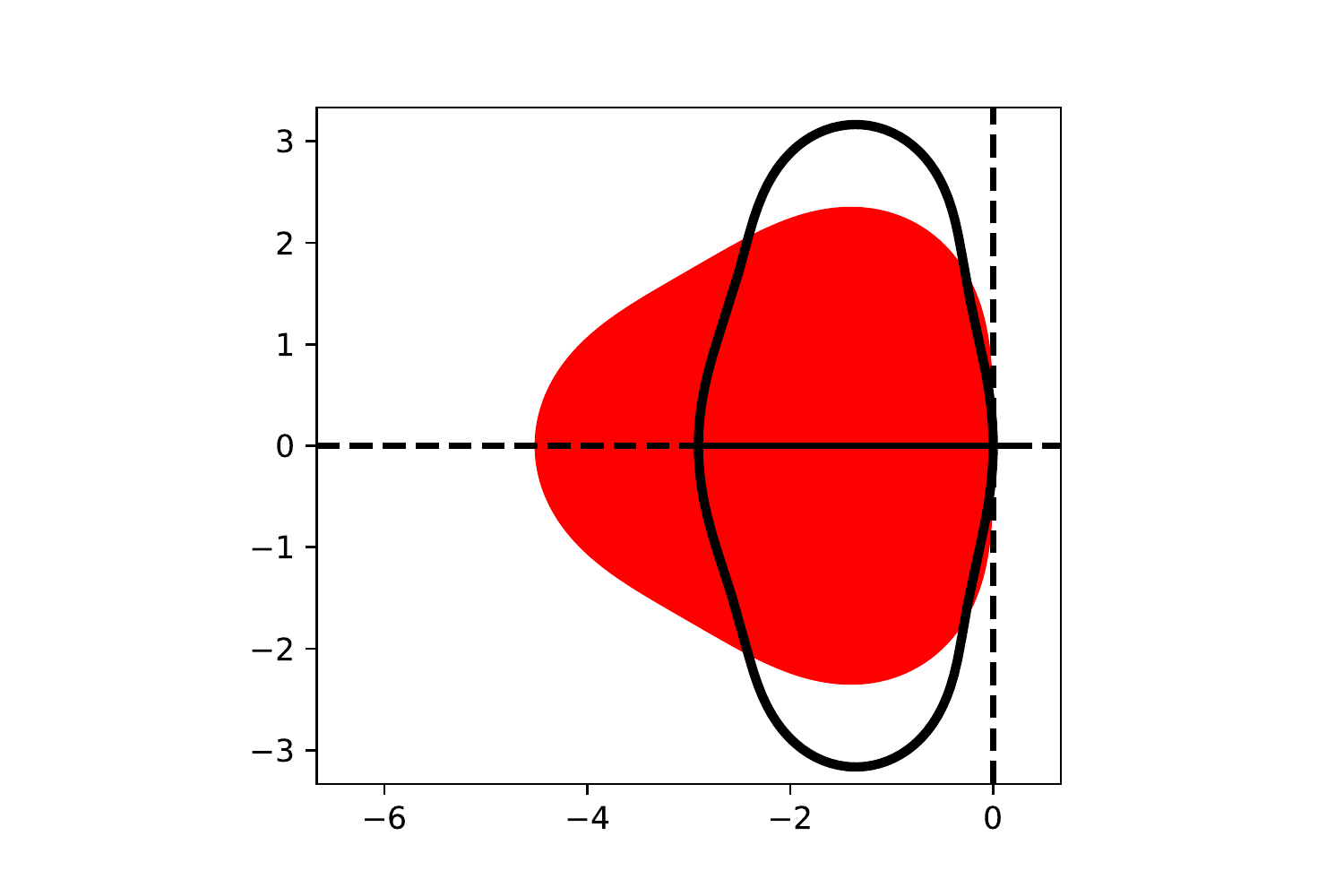}
        \caption{SSPERK$(3,2)-\wt$}
        \label{fig:SSPERK32}
    \end{subfigure}
    \caption{The stability regions~(red) of SSPERK$(2,2)$~(Figure-\ref{fig:SSPERK22BestNew01}) and SSPERK$(3,2)$~(Figure-\ref{fig:SSPERK32}) methods and the stability regions of the associated first order pair~(black contour).}\label{fig:NumericalStabRegSSPERK(s,2)}
\end{figure}

\subsection{Embedded pairs for SSPERK(s,3) methods}\label{sec:SSPERK(s,3)}
As a next step, we give embedded pairs for SSPERK$(s,3)$ methods. Kraaijevanger made the first step to determine the SSPERK$(s,3)$ methods. Namely, he treated the SSPERK$(4,3)$ case~\cite{K91}. Later the SSPERK$(s,3)$ methods were characterized, where $s=n^2$ and $n\geq 2$ is an integer~\cite{K08}. The corresponding SSP coefficients are $n^2-n=s-\sqrt{s}$. The Butcher form of SSPERK$(s,3)$ methods is given \eqref{eq:ButcherASSPERK(s,3)}-\eqref{eq:ButcherbSSPERK(s,3)}.
\begin{equation}\label{eq:ButcherASSPERK(s,3)}
    A = \begin{tikzpicture}[baseline={(m.center)}]
            \matrix [matrix of math nodes,left delimiter=(,right delimiter=)] (m)
            {0 &  &    &  &    &  &  &  \\
             \frac{1}{n(n-1)}  &  &    &  &    &  &  &  \\
             \frac{1}{n(n-1)}  &  \ddots  &    &  &    &  &  &  \\
                \frac{1}{n(n-1)} & &  \ddots&  &  &  &  &  \\
                 \vdots    &   & & \ddots &  &  &   &  &  &  \\
                \vdots       & &  & & \ddots  &  &  &  \\
                \frac{1}{n(n-1)} & \ldots & \frac{1}{n(n-1)}   & \ldots  & \ldots &  \frac{1}{n(n-1)}  &   &   &   & \\
                \frac{1}{n(n-1)}          &  & \frac{1}{n(n-1)} & \frac{1}{n(2n-1)}  & \ldots  & \frac{1}{n(2n-1)}  & \frac{1}{n(n-1)}   & &   &   \\
                \vdots    &   & \vdots & \vdots &   &  \vdots & \vdots & \ddots  \\
               \undermat{\frac{(n-2)(n-1)}{2}}{\frac{1}{n(n-1)} & \ldots &} \frac{1}{n(n-1)}  &\frac{1}{n(2n-1)}   &  \ldots & \frac{1}{n(2n-1)}  & \undermat{\frac{n(n-1)}{2}-1}{\frac{1}{n(n-1)} & \ldots &}  \frac{1}{n(n-1)}   & 0 \\
            };      
            \draw[loosely dotted] (m-1-1)-- (m-10-10);
            \draw (m-8-4.north west) rectangle (m-10-6.south east);
            \end{tikzpicture}\in\mathbb{R}^{n^2\times n^2},
\end{equation}

\noindent where the submatrix in the rectangle is a $\left(\frac{n(n-1)}{2}\right)\times\left(2n-1\right)$ dimensional matrix and
\begin{align}\label{eq:ButcherbSSPERK(s,3)}
b^T=\left[\underbrace{\frac{1}{n(n-1)},\ldots,\frac{1}{n(n-1)}}_{\frac{(n-1)(n-2)}{2}},\underbrace{\frac{1}{n(2n-1)},\ldots,\frac{1}{n(2n-1)}}_{2n-1},\underbrace{\frac{1}{n(n-1)},\ldots,\frac{1}{n(n-1)}}_{\frac{n(n-1)}{2}}\right]\in\mathbb{R}^{n^2}.
\end{align}

A possible embedded pair for the SSPERK$(4,3)$ method can be found in \cite{Ketcheson11} (Example 6.1.). The construction of this embedded pair $\bt_1=\left[\frac{1}{3},\frac{1}{3},\frac{1}{3},0\right]$ is related to the SSPERK$(3,2)$ method. Taking into account the Butcher forms of the SSPERK$(s,2)$ and SSPERK$(s,3)$ methods, one can realize that this kind of embedded pair can be only achieved in this exceptional case. Taking into account the desired properties (i)-(iv) our numerical searches suggest the embedded pair 
$$
\bt_2=\left[\frac{1}{4},\frac{1}{4},\frac{1}{4},\frac{1}{4}\right].
$$

Comparing the stability regions of the embedded pairs $\bt_1$ and $\bt_2$ it is easy to see that the embedded pair $\bt_2$ has a significantly larger stability region in Figure \ref{fig:StabRegSSPERK(4,3)}.

\begin{figure}[htbp]
    \centering
    \begin{subfigure}[b]{0.30\textwidth}
        \includegraphics[width=\textwidth,height=\textwidth]{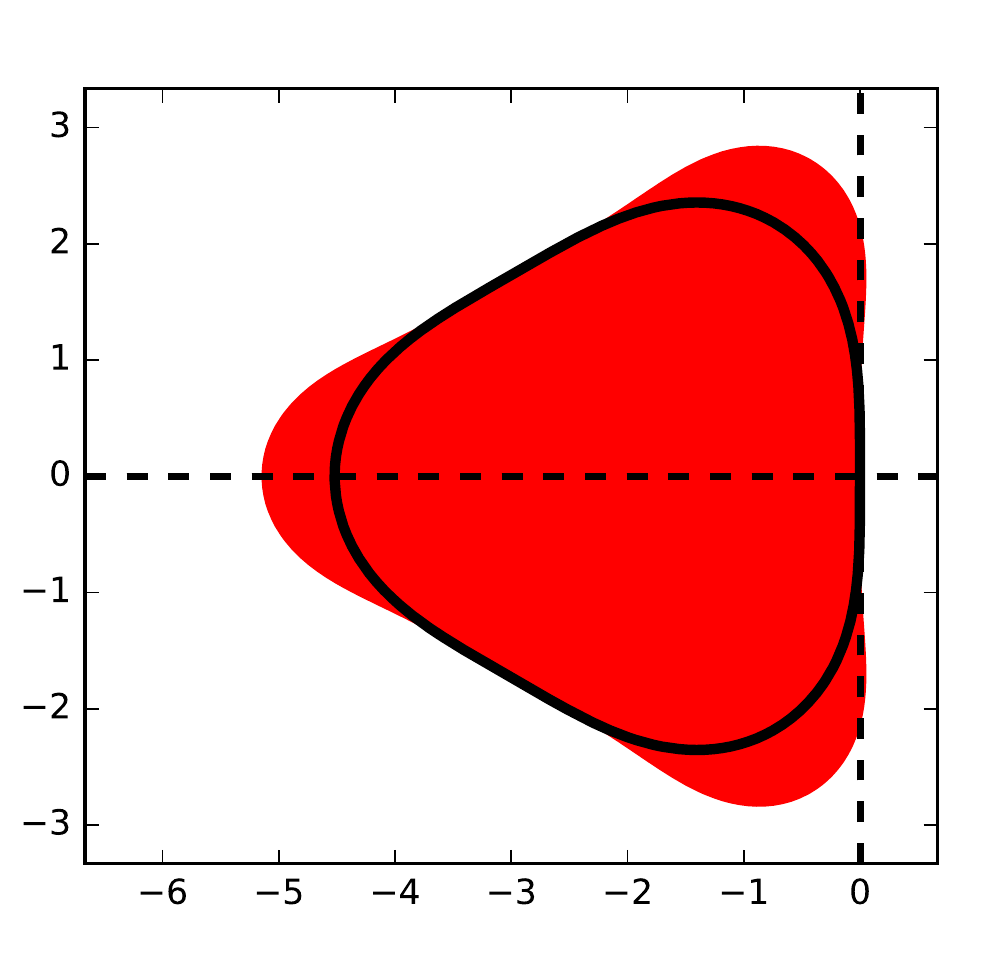} \caption{SSPERK$(4,3)-\bt_1$}
        \label{fig:SSPERK43b}
    \end{subfigure}~ 
    \begin{subfigure}[b]{0.30\textwidth}
        \includegraphics[width=\textwidth,height=\textwidth]{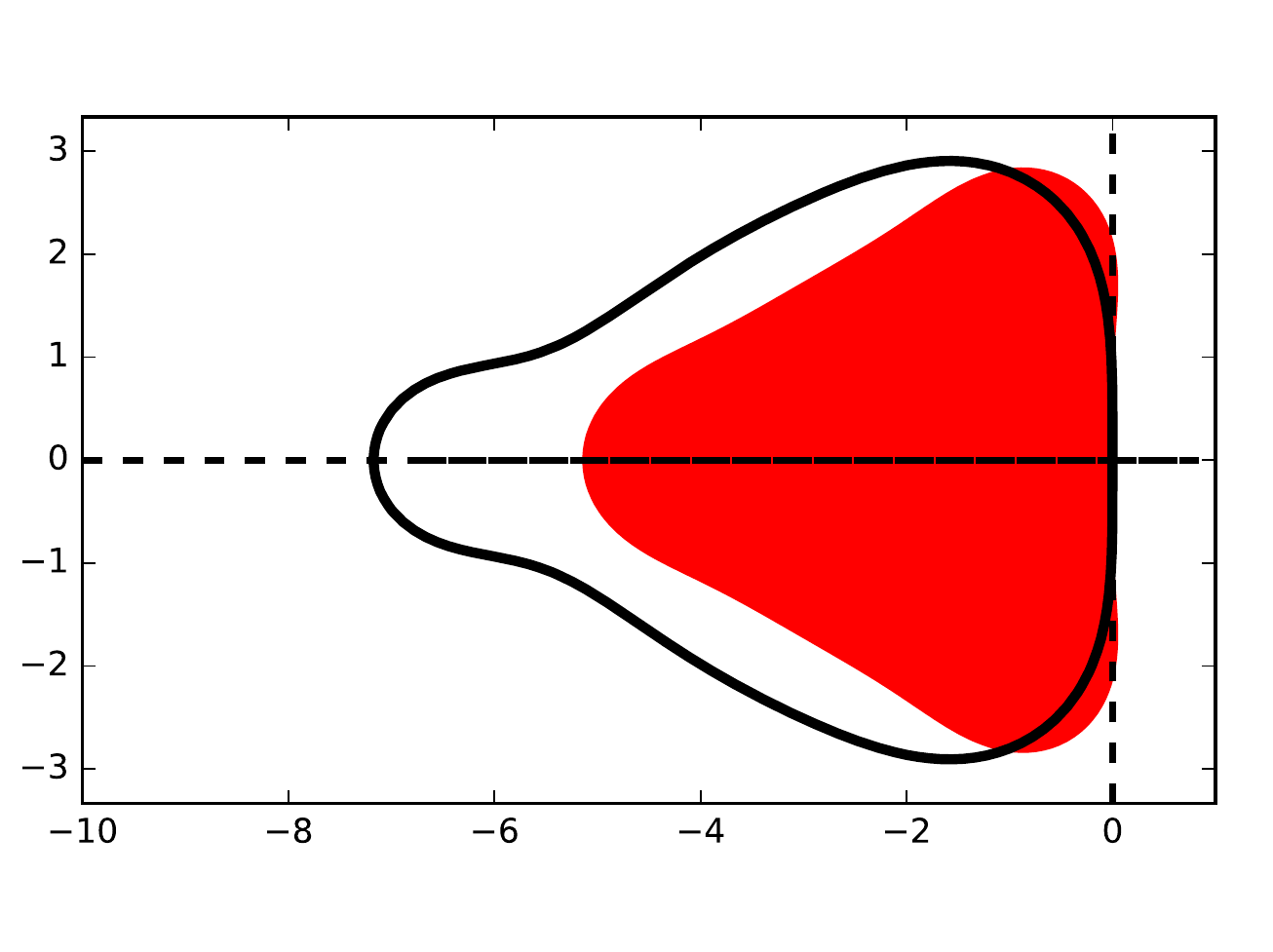}
        \caption{SSPERK$(4,3)-\bt_2$}
       \label{fig:SSPERK43}
    \end{subfigure}
\caption{The stability regions of SSPERK$(4,3)$ method (red) and the black contours of the embedded SSPERK$(4,2)$ methods.}\label{fig:StabRegSSPERK(4,3)}
\end{figure}

\begin{figure}[htbp]
    \centering
    \begin{subfigure}[b]{0.47\textwidth}
        \includegraphics[width=\textwidth,height=\textwidth]{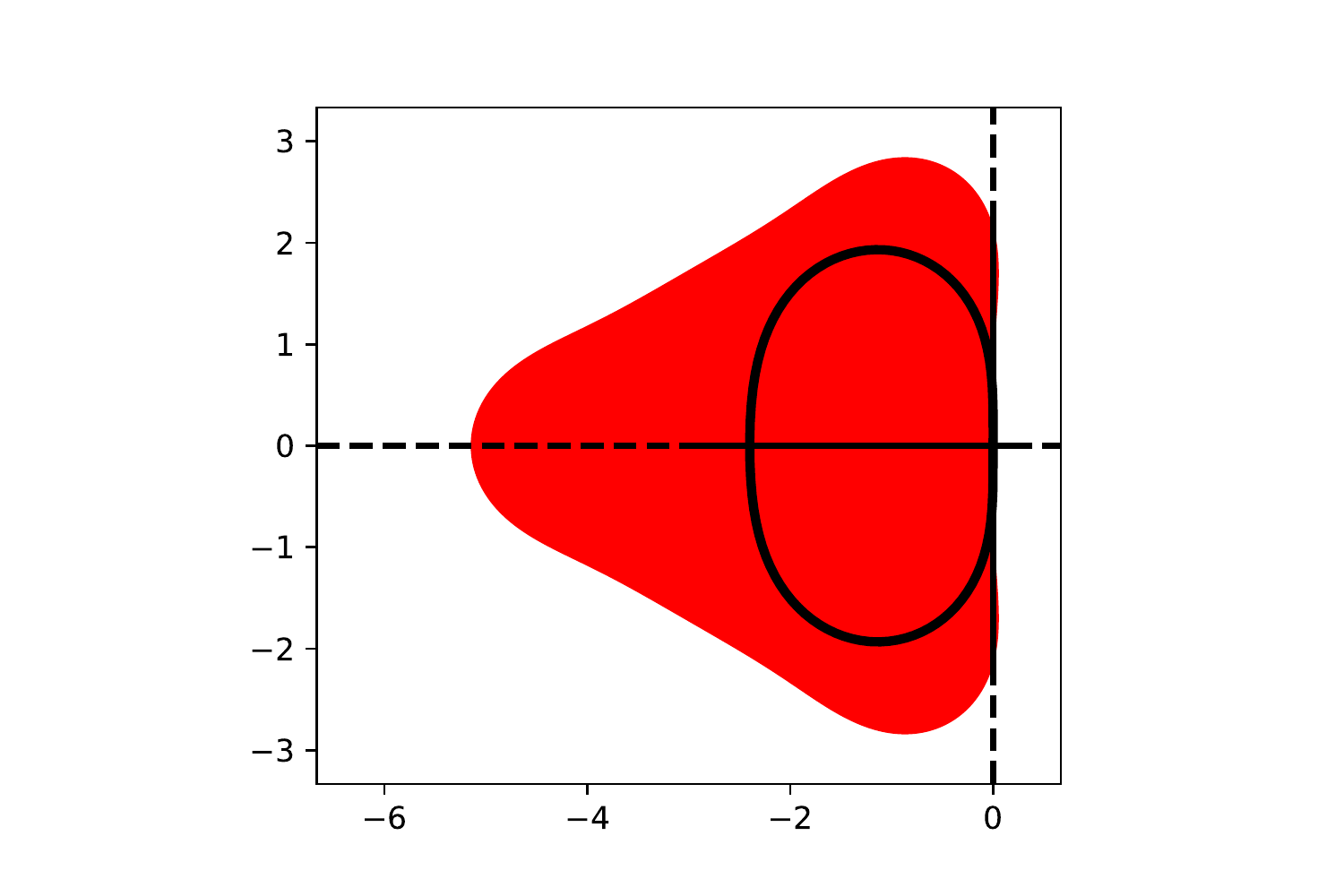} \caption{SSPERK$(4,3)-\wt$}
        \label{fig:SSPERK43Best}
    \end{subfigure}
    \begin{subfigure}[b]{0.47\textwidth}
        \includegraphics[width=\textwidth,height=\textwidth]{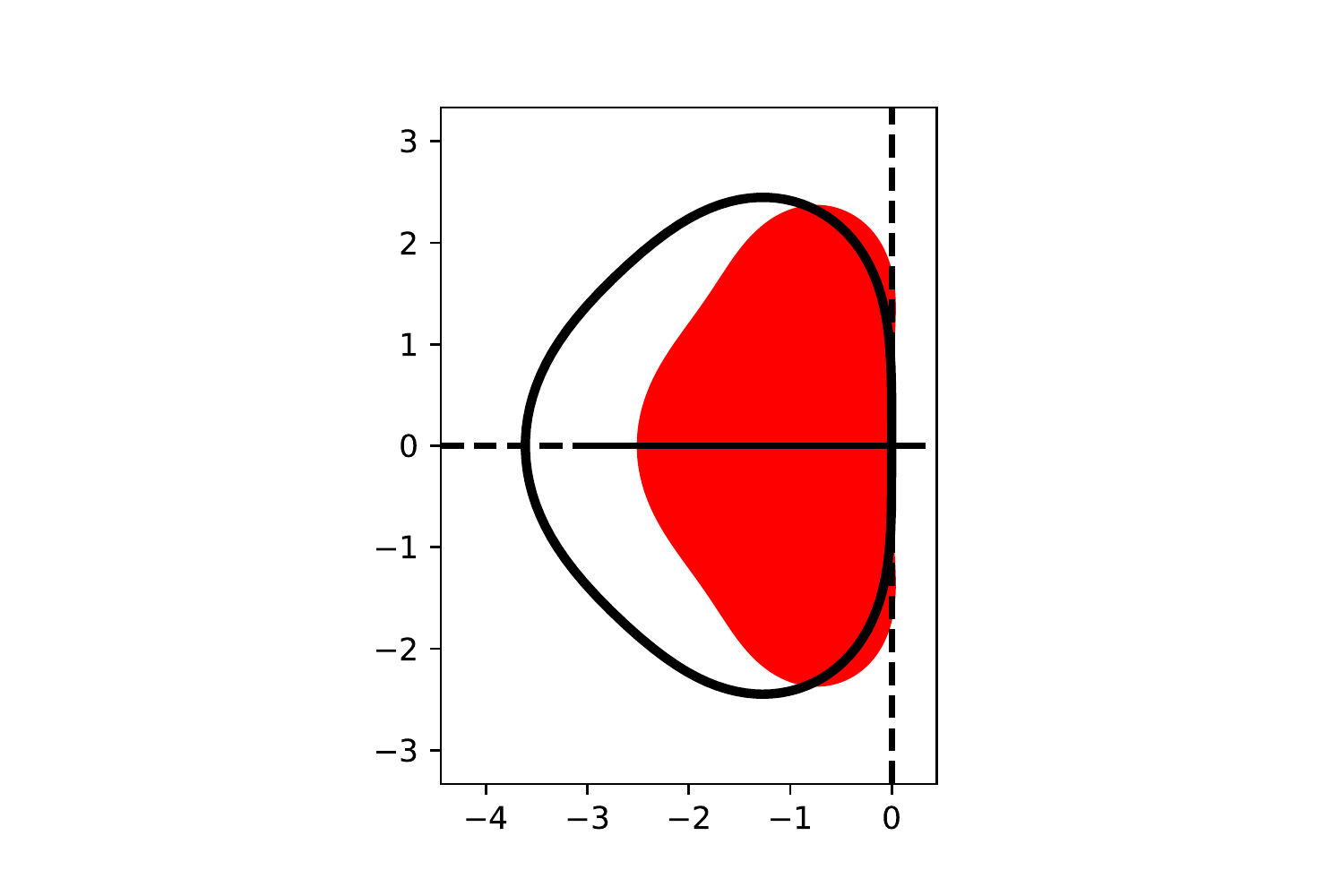}
        \caption{SSPERK$(3,3)-\wt$}
        \label{fig:SSPERK33}
    \end{subfigure}
    \caption{The stability regions of SSPERK$(s,3)$~(Figure-\ref{fig:SSPERK43Best}) and SSPERK$(s,3)$~(Figure-\ref{fig:SSPERK33}) methods (red) and the black contours of the embedded ERK$(s,1)$ methods with $s = 4, 3$ respectively~(obtained from the numerical optimization~\ref{eq:numOptimization}).}\label{fig:NumericalStabRegSSPERK(s,3)}
\end{figure}
Now, we turn our attention to SSPERK$(n^2,3)$ methods, where $n\geq 3$ is an integer. Numerical searches failed to find any embedded pairs with $\tilde{\sspC}=n^2-n$. Therefore, taking into account a general simple structure (iii) and the desired property (iv), we are simply looking for non-defective embedded ERK$(n^2,2)$ methods. Similarly to the SSPERK$(4,3)$ case, the suggested non-defective embedded pair is
$$
\tilde{b}^T=\left[\frac{1}{n^2},\ldots,\frac{1}{n^2}\right].
$$

To demonstrate the efficiency of this choice we plot the stability regions of the SSPERK$(n^2,3)$ methods and the embedded method for different cases in Figure \ref{fig:StabRegSSPERK(s,3)}. Furthermore, we also give the stability radius measurement and error measurement values in~\cite{CondeFekete18}. The embedded pair is simply denoted by $\bt$. Our computations also show that the suggested embedded method has $\tilde{\sspC}>0$. In addition, it is optimal for the SSPERK$(4,3)$ case. The corresponding $R(\psi)$ values are also given in~\cite{CondeFekete18}.

\begin{figure}[ht!]
    \centering
    \begin{subfigure}[b]{0.25\textwidth}
        \includegraphics[width=\textwidth,height=\textwidth]{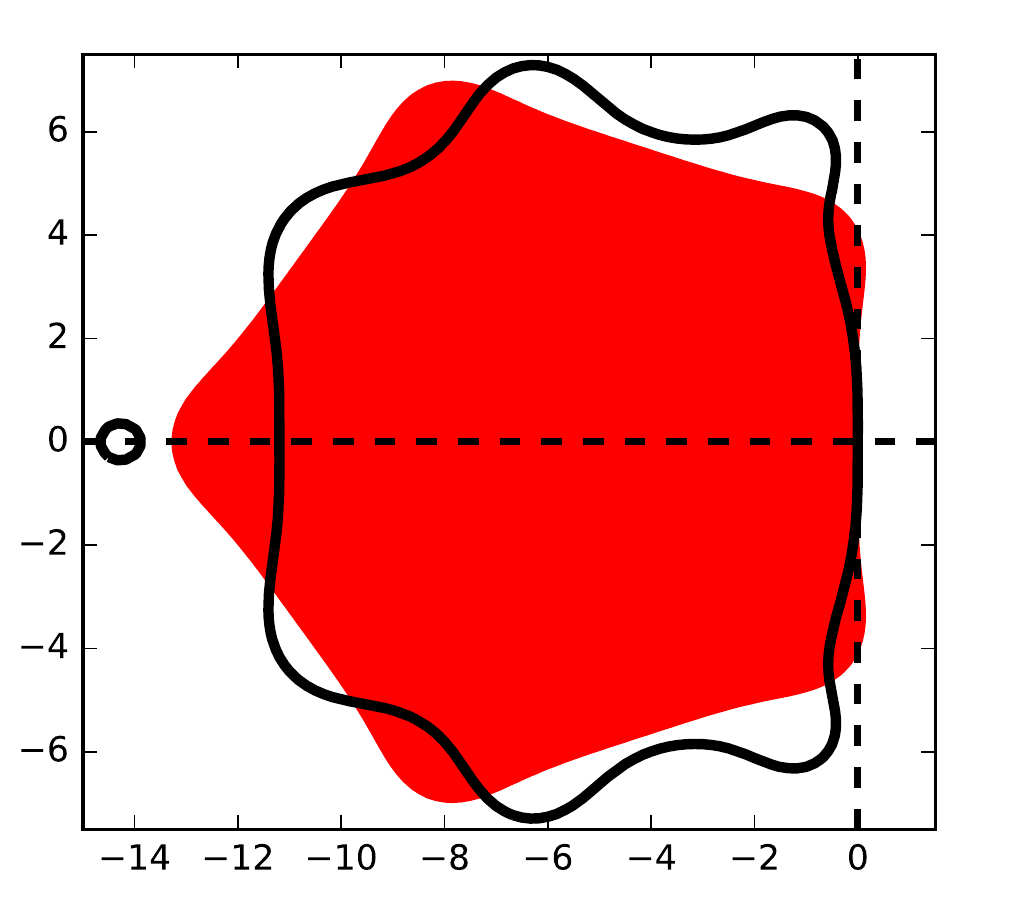} \caption{SSPERK$(9,3)$}
        \label{fig:SSPERK93}
    \end{subfigure}~ 
    \begin{subfigure}[b]{0.25\textwidth}
        \includegraphics[width=\textwidth,height=\textwidth]{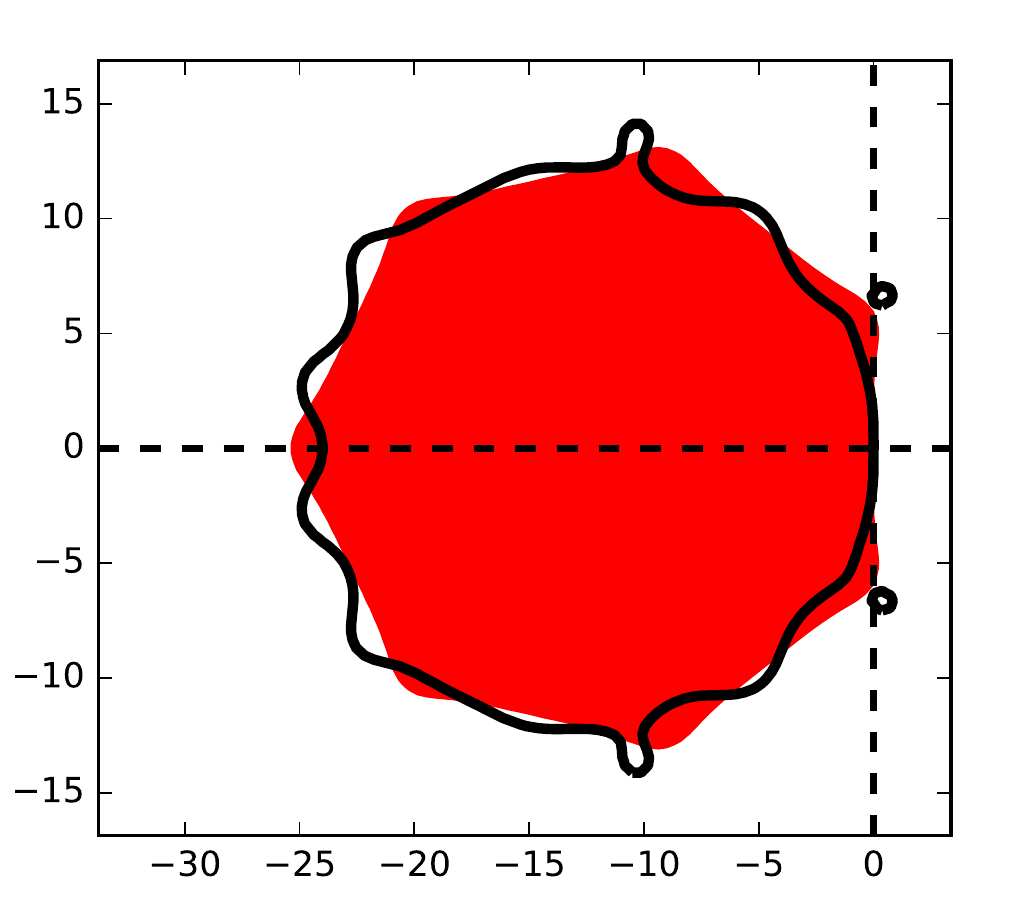} \caption{SSPERK$(16,3)$}
        \label{fig:SSPERK163}
    \end{subfigure}~
        \begin{subfigure}[b]{0.25\textwidth}
        \includegraphics[width=\textwidth,height=\textwidth]{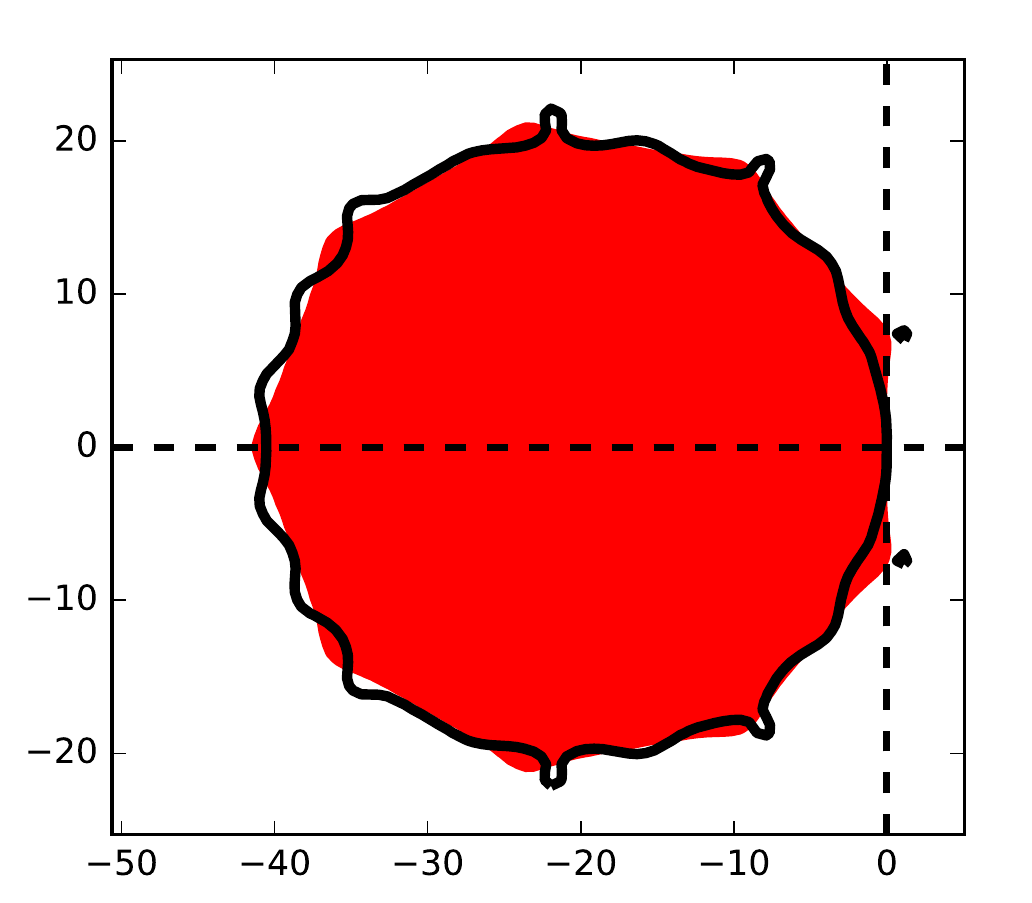} \caption{SSPERK$(25,3)$}
        \label{fig:SSPERK253}
    \end{subfigure}~ 
    \begin{subfigure}[b]{0.25\textwidth}
        \includegraphics[width=\textwidth,height=\textwidth]{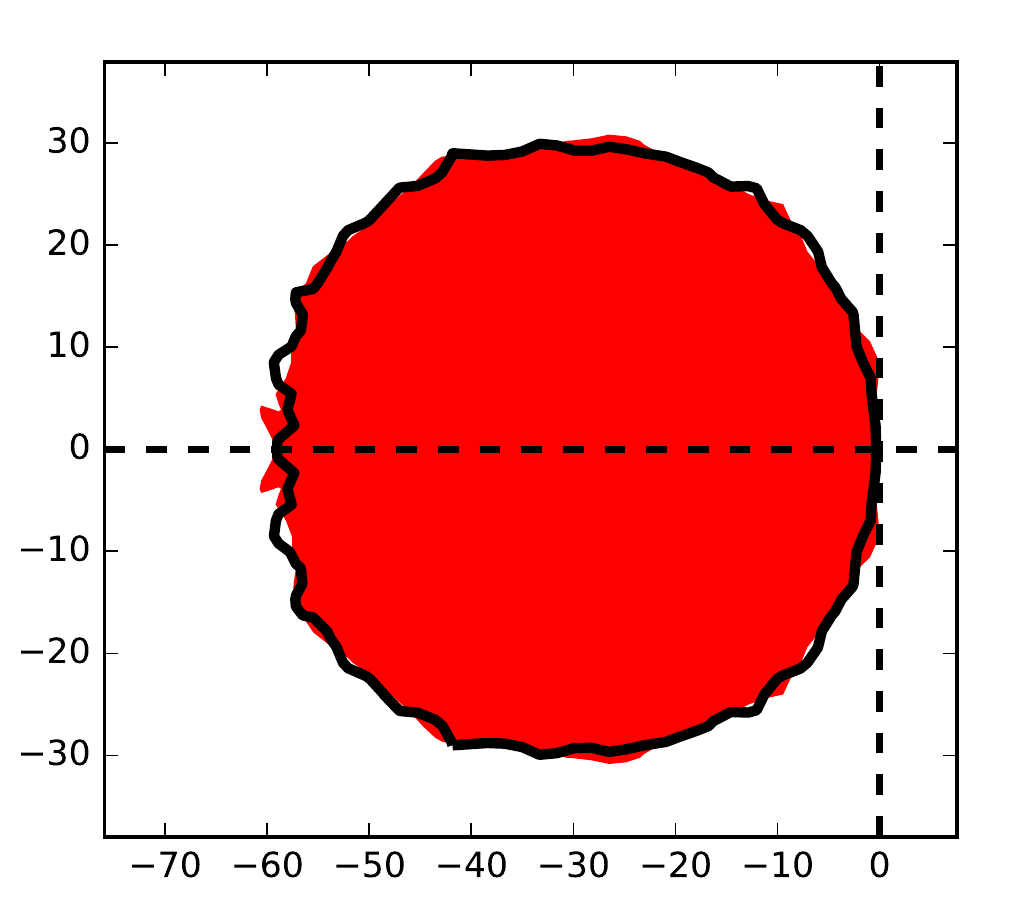} \caption{SSPERK$(36,3)$}
        \label{fig:SSPERK363}
    \end{subfigure}
    \caption{The stability regions of the SSPERK$(n^2,3)$ methods (red) and the contours of the embedded SSPERK$(n^2,2)$ method with $\bt=\frac{1}{n^2}\mathbf{e}$ (black) for $n^2=9,16,25, 36$.}\label{fig:StabRegSSPERK(s,3)}
\end{figure}



\subsection{Embedded pairs for the SSPERK(10,4) method}\label{sec:SSPERK(10,4)}
In this section we consider the SSPERK$(10,4)$ method which has simple rational coefficients~\cite{Ketcheson11}. This method is popular in the literature and applications because it allows for a low-storage implementation~\cite{K08}. Furthermore, it has $\sspC =6$. The corresponding Butcher matrix of SSPERK$(10,4)$ is
\begin{equation*}
    A = \begin{tikzpicture}[baseline={(m.center)}]
            \matrix [matrix of math nodes,left delimiter=(,right delimiter=)] (m)
            {0 &  &    &  &    &  &  &  \\
             \frac{1}{6}  &  &    &  &    &  &  &  \\
              \frac{1}{6}  & \frac{1}{6}  &    &  &    &  &  &  \\
               \frac{1}{6}  & \frac{1}{6} & \frac{1}{6} &  &  &  &  &  \\
               \frac{1}{6}   &  \frac{1}{6}  & \frac{1}{6} &\frac{1}{6}  &  &  &   &  &  &  \\
              \frac{1}{15}       & \frac{1}{15}  & \frac{1}{15}   & \frac{1}{15}  & \frac{1}{15}   &  &  &  \\
 \frac{1}{15}       & \frac{1}{15}  & \frac{1}{15}   & \frac{1}{15}  & \frac{1}{15}   & \frac{1}{6} &  &  \\
  \frac{1}{15}       & \frac{1}{15}  & \frac{1}{15}   & \frac{1}{15}  & \frac{1}{15}   & \frac{1}{6} &   \frac{1}{6} &  \\
   \frac{1}{15}       & \frac{1}{15}  & \frac{1}{15}   & \frac{1}{15}  & \frac{1}{15}   & \frac{1}{6} &   \frac{1}{6} &  \frac{1}{6}\\
               \frac{1}{15} &   \frac{1}{15} & \frac{1}{15}&  \frac{1}{15}&  \frac{1}{15} &   \frac{1}{6}  & \frac{1}{6}  &        \frac{1}{6}& \frac{1}{6} & 0 \\
            };      
            \draw[loosely dotted] (m-1-1)-- (m-10-10);
            \end{tikzpicture}
\end{equation*}
and its Butcher array is
       $$
b^T=\left[\frac{1}{10},\frac{1}{10},\frac{1}{10},\frac{1}{10},\frac{1}{10},\frac{1}{10},\frac{1}{10},\frac{1}{10},\frac{1}{10},\frac{1}{10}\right].
$$


It turns out that we cannot expect non-defective embedded methods for the SSPERK(10,4) method since condition \eqref{eq:ordercond4c} is always satisfied, thus we want to find embedded ERK$(10,3)$ methods which violate the other three fourth-order conditions. However, our numerical searches failed to find embedded pairs which has $\sspC=6$. Therefore, we are looking for embedded ERK$(10,3)$ methods. During our investigation we have found eight potential embedded pairs which have nice structures. 
The potential embedded pairs are:
\begin{align*}
\bt_1&=\left[0, \frac{3}{8}, 0, \frac{1}{8}, 0, 0, 0, \frac{3}{8}, 0, \frac{1}{8}\right] &
\bt_2&=\left[\frac{3}{14}, 0, 0, \frac{2}{7}, 0, 0, 0, \frac{3}{7}, 0, \frac{1}{14}\right]\\
\bt_3&=\left[0, \frac{2}{9}, 0, 0, \frac{5}{18}, \frac{1}{3}, 0, 0, 0, \frac{1}{6}\right] &
\bt_4&=\left[\frac{1}{5}, 0, 0, \frac{3}{10}, 0, 0, \frac{1}{5}, 0, \frac{3}{10}, 0\right]\\
\bt_5&=\left[\frac{1}{10}, 0, 0, \frac{2}{5}, 0,  \frac{3}{10}, 0, 0, 0, \frac{1}{5}\right] &
\bt_6&=\left[\frac{1}{6}, 0, 0, 0, \frac{1}{3}, \frac{5}{18}, 0, 0, \frac{2}{9}, 0\right]\\
\bt_7&=\left[0, \frac{2}{5}, 0, \frac{1}{10}, 0, 0, 0, \frac{1}{5}, \frac{3}{10}, 0\right] &
\bt_8&=\left[\frac{1}{7}, 0, \frac{5}{14}, 0, 0, 0, 0, \frac{3}{14}, \frac{2}{7}, 0\right]\\
\end{align*}
We would like to note that we have good candidates in this case but not an outstanding one. Taking into account their stability radius and error control measurement values~\cite{CondeFekete18}, we suggest the embedded pair $\bt_3$.
In order to demonstrate our choice we plot the stability regions of the embedded pairs in Figure \ref{fig:StabRegSSPERK(10,4)}.

\begin{figure}[htbp]
\centering
        \begin{subfigure}[b]{0.20\textwidth}
        \includegraphics[width=\textwidth]{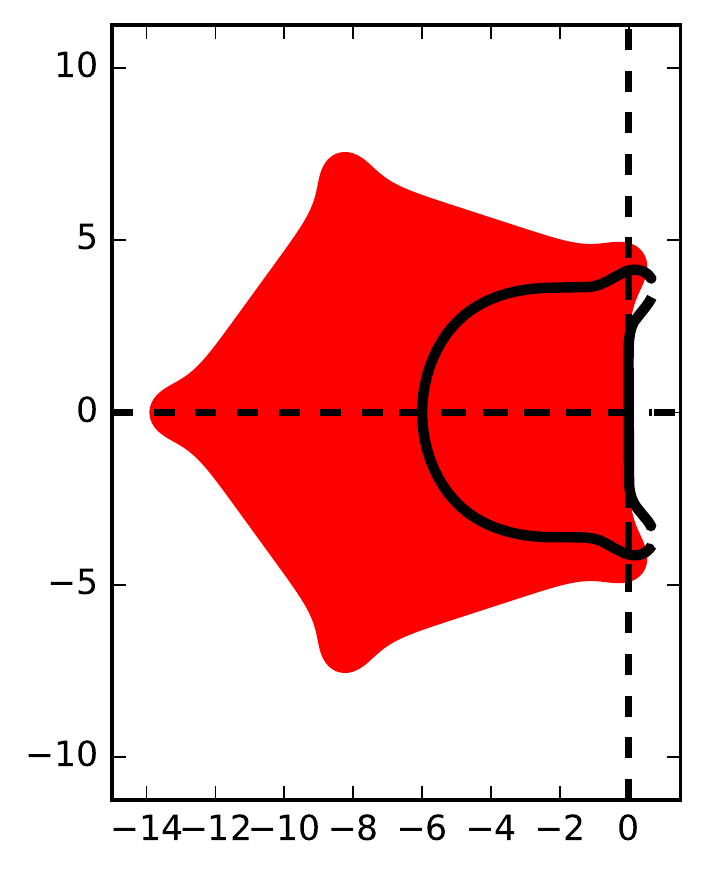}
        \caption{SSPERK$(10,4)~\bt_1$}
        \label{fig:SSPERK104b1}
    \end{subfigure}~ 
        \begin{subfigure}[b]{0.20\textwidth}
        \includegraphics[width=\textwidth]{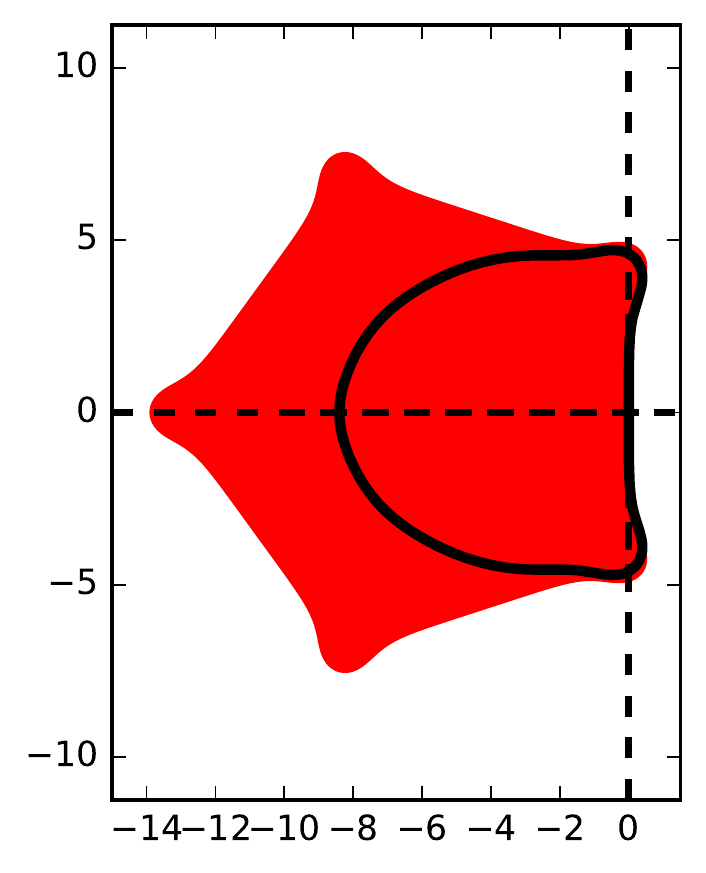}
        \caption{SSPERK$(10,4)~\bt_2$}
        \label{fig:SSPERK104b2}
    \end{subfigure}~
        \begin{subfigure}[b]{0.20\textwidth}
        \includegraphics[width=\textwidth]{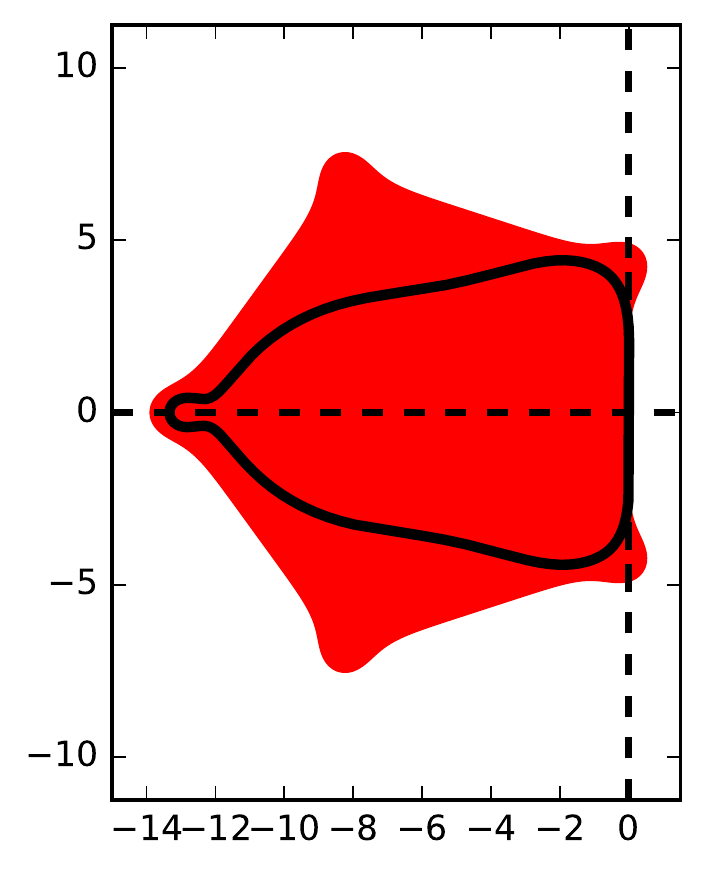}
         \caption{SSPERK$(10,4)~\bt_3$}
        \label{fig:SSPERK104b3}
    \end{subfigure}~
        \begin{subfigure}[b]{0.20\textwidth}
        \includegraphics[width=\textwidth]{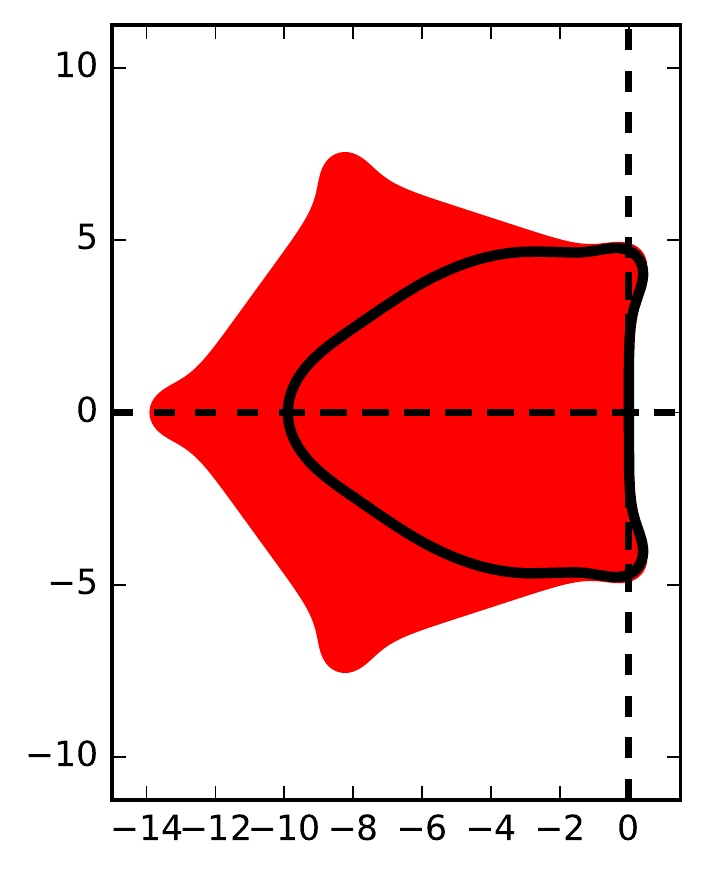}
        \caption{SSPERK$(10,4)~\bt_4$}
        \label{fig:SSPERK104b4}
    \end{subfigure}\\
        \begin{subfigure}[b]{0.20\textwidth}
        \includegraphics[width=\textwidth]{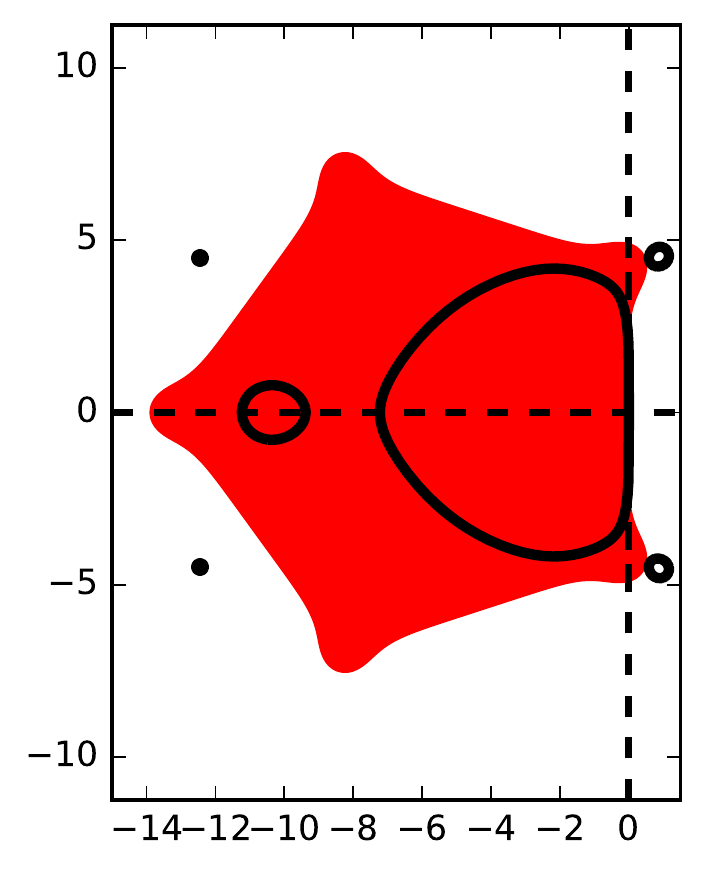}
        \caption{SSPERK$(10,4)~\bt_5$}
        \label{fig:SSPERK104b5}
    \end{subfigure}~ 
        \begin{subfigure}[b]{0.20\textwidth}
        \includegraphics[width=\textwidth]{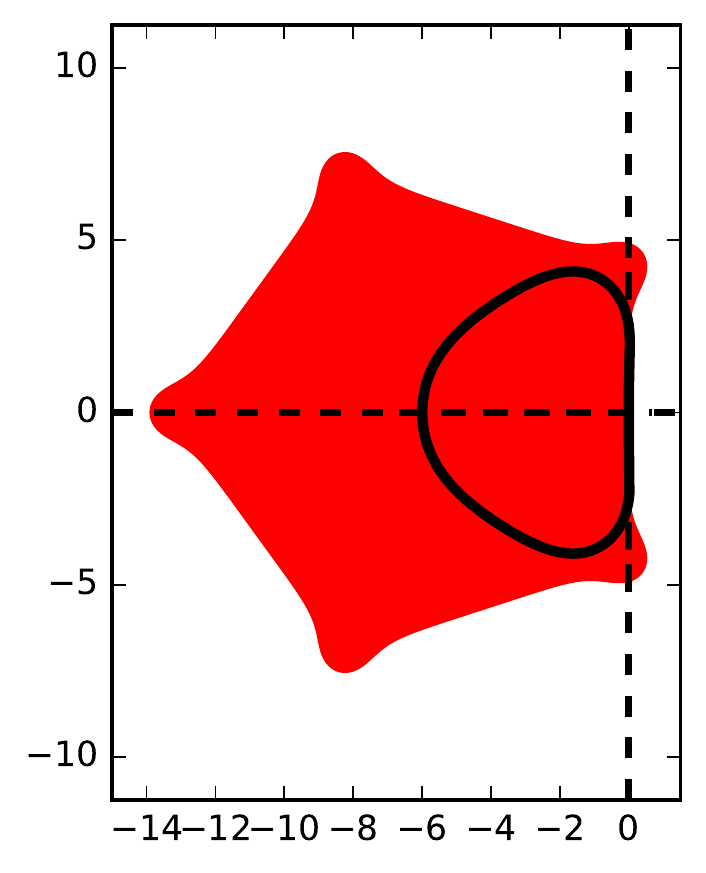}
        \caption{SSPERK$(10,4)~\bt_6$}
        \label{fig:SSPERK104b6}
    \end{subfigure}~
        \begin{subfigure}[b]{0.20\textwidth}
        \includegraphics[width=\textwidth]{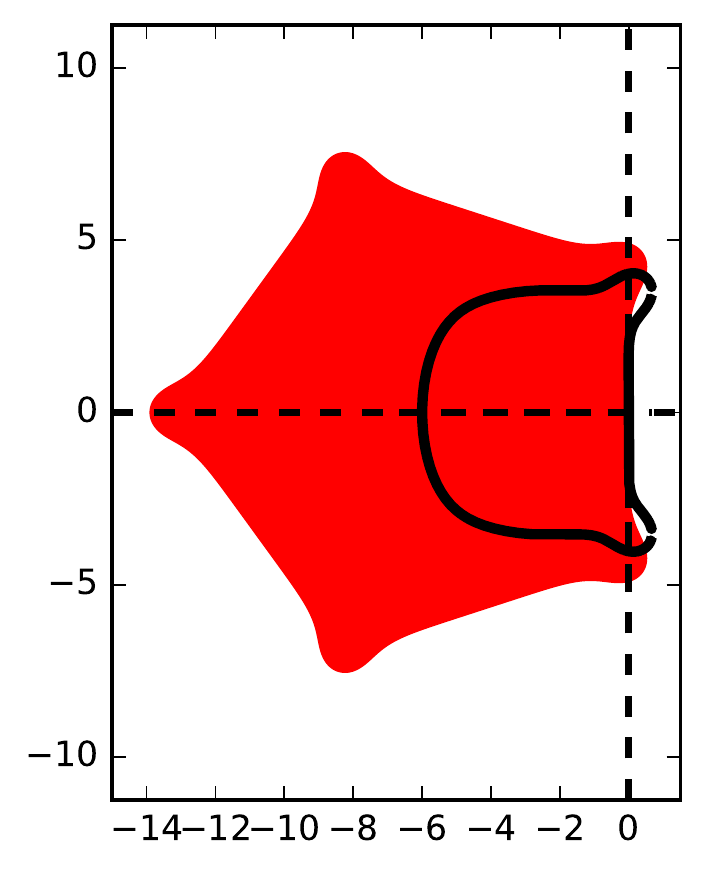}
         \caption{SSPERK$(10,4)~\bt_7$}
        \label{fig:SSPERK104b7}
    \end{subfigure}~
        \begin{subfigure}[b]{0.20\textwidth}
        \includegraphics[width=\textwidth]{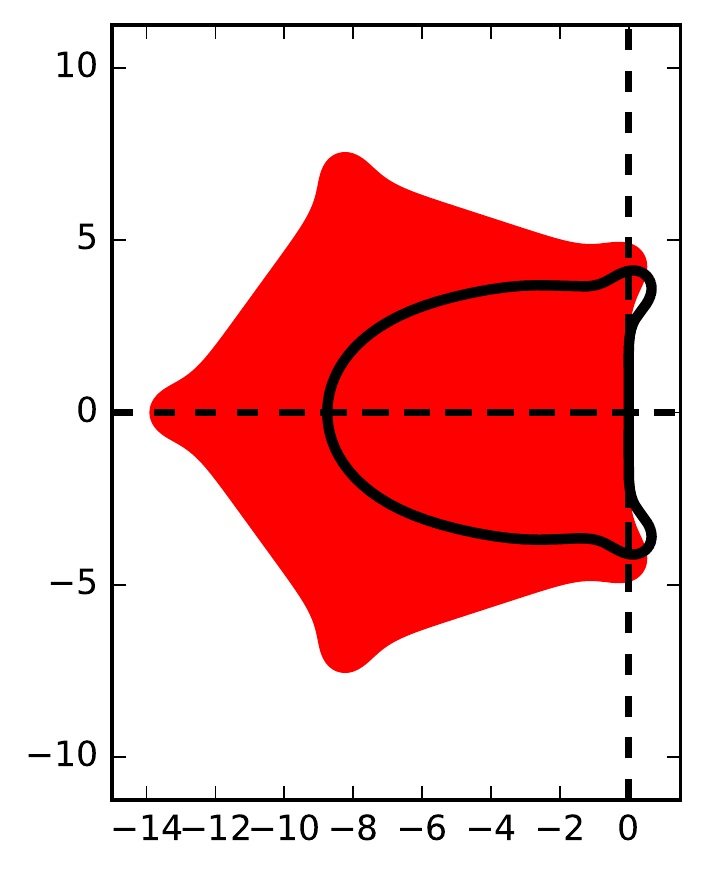}
        \caption{SSPERK$(10,4)~\bt_8$}
        \label{fig:SSPERK104b8}
    \end{subfigure}

    \caption{The stability regions of the SSPERK$(10,4)$ method (red) and the black contours of the embedded SSPERK$(10,3)$ methods.}\label{fig:StabRegSSPERK(10,4)}
\end{figure}

\begin{figure}[htbp!]
    \centering
    \begin{subfigure}[b]{0.5\textwidth}
        \includegraphics[width=\textwidth,height=\textwidth]{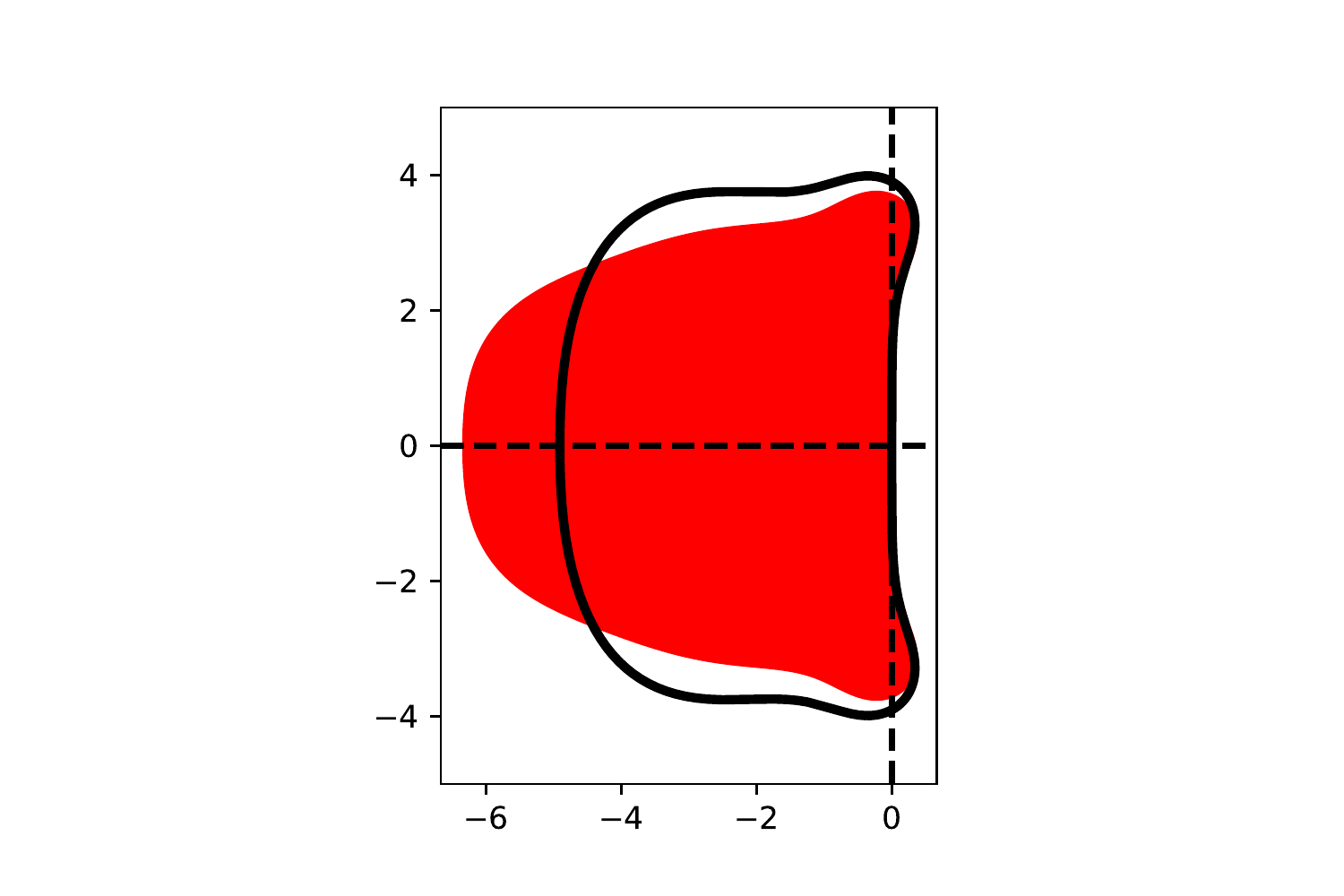} \caption{SSPERK$(6,4)-\wt$}
        \label{fig:SSPERK64BestNew01}
    \end{subfigure}
    \caption{The stability regions of SSPERK$(6,4)$ methods (red) and the black contour of the embedded ERK$(6,3)$ method~(obtained from the numerical optimization~\ref{eq:numOptimization}).}\label{fig:NumericalStabRegSSPERK(6,4)}
\end{figure}

\section{Numerical Results}\label{sec:numericalresults}


\subsection{Step-size control}


Whenever a starting step-size $\h$ has been chosen,
the adaptive method computes two approximations to the solution, $u_{n+1}$ and $\hat{u}_{n+1}$. Then an estimate of the local error for the less precise result is $u_{n+1} - \hat{u}_{n+1}$ \cite{HairerNorsettWanner93,Ascher98}.
Adaptive time-stepping techniques aim for this error to satisfy componentwise
\begin{equation}
| u_{n+1} - \hat{u}_{n+1} | \leq sc, \qquad sc = Atol + \max\left(|u_n|, | u_{n+1} | \right)\cdot Rtol
\label{eq:error_sc}
\end{equation}

\noindent where $Atol$ and $Rtol$ defines respectively the desired absolute and relative tolerance. The estimated local truncation error is thus defined as: 

$$ \err{n+1} = \left\| \frac{(u_{n+1} - \hat{u}_{n+1})}{sc} \right\|_\infty $$

\noindent where $sc$ is given by~\eqref{eq:error_sc}. The error is inspected and a decision is made whether to accept the computed solution or not. The optimal step-size is obtained by

\begin{equation}
\h_{opt} = \h \cdot \beta_{n+1} 
\label{eq:optionalStepSize}
\end{equation}

\noindent where $\beta_{n+1}$ is determined by the choice of error control algorithm. Often the optimal step-size $\h_{opt}$ is scaled conservatively by a safety factor. We chose a safety factor of $fac = 0.9, facmin = 0.1$. Furthermore, $\h$ is neither allowed to increase nor to decrease too fast. So following~\cite{HairerNorsettWanner93}

\begin{equation}
    \h_{opt} = \h \cdot \min( facmax, \max(facmin, fac \cdot \beta_{n+1}))
\label{eq:optionalStepSize1}
\end{equation}

If $\err{n+1} \leq 1$, then the computed solution $\y_{n+1}$ is accepted and advanced. The next time step is computed with $\h_{opt}$ as step-size. If $\err{n+1} > 1$, then step is rejected and the computations are repeated with the new step-size $\h_{opt}$ . We chose a maximal step-size increase factor of $facmax = 5$ to limit large increases in the step-size. In computing the approximate solution immidiately following a step-rejection, we set $facmax = fac = 0.9$ to prevent an infinite loop that we observed on rare occasions - similar strategies are employed by~\cite{ShampineReichelt97,bogackiShampine89,HairerNorsettWanner93}. 


\subsubsection*{Starting Step-Size} \label{sec:startingStepSize}
We follow the algorithm of Gladwell, Shampine \& Brankin~\cite{gladwell1987,HairerNorsettWanner93} which is based on the hypothesis that
$$\text{local error } \approx C(\Delta t)^{p+1}u^{(p+1)}(t_0).$$
The starting step-size is computed following the algorithm in~\cite[p. 169.]{HairerNorsettWanner93}. In the case of simple ODEs (Van der Pol and Brusselator), we used the computed \emph{optimal} starting step-size. For the PDEs, we take the minimum of the starting step-size returned from the algorithm and the step-size restriction based on the CFL constraints to avoid any stability issues from influencing the numerical study.

\subsubsection*{Asymptotic error control} \label{sec:error_control}

There exists a variety of error control algorithms~\cite{Kennedy2003ARS,soderlind1998,Soderlind2002,soderlind2003digital,soderlind2006time}. In the numerical experiment that follows, we used four separate error control strategies. The standard time adaptivity {I controller} provides a prospective time step estimate entirely based on the current local error estimate

\begin{equation}
    \beta_{n+1} = \err{n+1}^{-k_1 / p}.
\label{eq:IController}
\end{equation}

\noindent By default, we take $k_1 = 1$. The {PI controller} uses the two most recent local truncation errors estimated in its adaptivity algorithm:

\begin{equation}
    \beta_{n+1} = \err{n+1}^{-k_1 / p} \err{n}^{k_2 / p}.
 \label{eq:PIController}
 \end{equation}
 
\noindent Here, the default values are $k_1 = 0.8$, $k_2 = 0.31$. The {PID controller} uses the information from the three most recent time steps to provide an optimal step-size:

\begin{equation}
    \beta_{n+1} = \err{n+1}^{-k_1 / p} \err{n}^{k_2 / p} \err{n-1} ^{-k_3 / p}.
 \label{eq:pidController}
 \end{equation}
Here, the default values are $k_1 = 0.58$, $k_2 = 0.21$, $k_3 = 0.1$.

Lastly, {Explicit Gustafsson Controller} is primarily useful in combination with explicit Runge--Kutta methods and was proposed in \cite{gustafsson91},
 
 \begin{equation}
\beta =\displaystyle
\begin{cases}
 \err{1}^{-1/ p} & \text{ on the first step},\\
    \err{n+1}^{-k_1 / p} \left( \err{n+1} / \err{n} \right)^{k_2 / p} & \text{ on the subsequent steps}.
\end{cases}
 \label{eq:explicitGustafssonController}
 \end{equation}
 
 The values of $k_1 = 0.367$, $k_2 = 0.268$. The default values chosen are similarly used by SUNDIALS/CVODE~\cite{hindmarsh2005sundials}. In this estimate, a floor of $\err{n}> 1e -10$ is enforced to avoid division by zero errors.


In this section we present some experiments to test the numerical efficiency of the new embedded pairs constructed above when applied to several oscillatory IVPs and hyperbolic PDEs. We run the Runge--Kutta pairs on the test problems described below and for the range of tolerances $\Rtol = \Atol = 1e-2, 1e-3, \ldots, 1e-7$.


We compare several classical embedded pairs with the newly developped SSPERK embedded pairs. From the literature, we have the following methods:
\emph{RKF23}, \emph{Ceschino24}, \emph{RKF23b}, \emph{Fehlberg12}, \emph{Fehlberg12b}, \emph{Bogackishampine}, \emph{Merson45}, \emph{Zonneveld43}, \emph{Fehlberg45}, and \emph{Dormandprince54}.
The methods \emph{RKF23b}~\cite[Table 4.4b]{HairerNorsettWanner93} and \emph{RKF23}~\cite[Table 4.4a]{HairerNorsettWanner93}: the solution $\y$ advanced is second order with a third order embedded solution $\hat{\y}$ used for error estimation.
\emph{Ceschino24}~\cite[Table 4.1]{HairerNorsettWanner93}: $\y$ is of second order while the embedded solution is fourth order accurate.
\emph{Fehlberg12}~\cite[Table 4.7a]{HairerNorsettWanner93}: the solution advanced is first order accurate, while the embedded weight used for error approximation is of second order;
\emph{Fehlberg12b}~\cite[Table 4.7b]{HairerNorsettWanner93}: the embedded weight is second-order accurate, the method advanced numerically gives a slightly higher order than 1;
\emph{Bogackishampine}~\cite{bogackiShampine89} is of order three with four stages with the First Same As Last (FSAL) property, so that it uses approximately three function evaluations per step (although we do not exploit this property in this study). This method is implemented in the \emph{ODE23} function in MATLAB.
\emph{Merson45}~\cite[Table 4.1]{HairerNorsettWanner93}: often implemented in NAG softwares, $\hat{\y}_1$ is of order 5 for linear equations with constant coefficients; for nonlinear problems it is of order 3.
\emph{Zonneveld43}~\cite[Table 4.2]{HairerNorsettWanner93}: fourth-order accurate weight used to advance the solution while a third-order embedded weight is used for error estimation.
\emph{Fehlberg45}~\cite{HairerNorsettWanner93} is a method of order four with an error estimator of order five.
\emph{Dormandprince54}~\cite{HairerNorsettWanner93,dormandPrince198167} uses six function evaluations to calculate fourth- and fifth-order accurate solutions and is currently the default method in \emph{ODE45} solver.

We report the error coefficients for the all the second; third; and fourth-order explicit embedded pairs in~\cite{CondeFekete18}. The global errors are calculated by using a very accurate solution calculated by  MATLAB's \emph{ODE45} routine with tolerances set to {\tt AbsTol= }{\tt  RelTol= } $10^{-14}$. The numerical tests were conducted by~\cite{CondeFekete18} and the numerical tests for the oscillatory IVPs were reproduced in~\cite{tempus17}.


\subsubsection{Step-size Control: ODE Tests}

\begin{figure}[h]
    \centering
    \begin{subfigure}[b]{0.475\textwidth}
        \includegraphics[width=\textwidth]{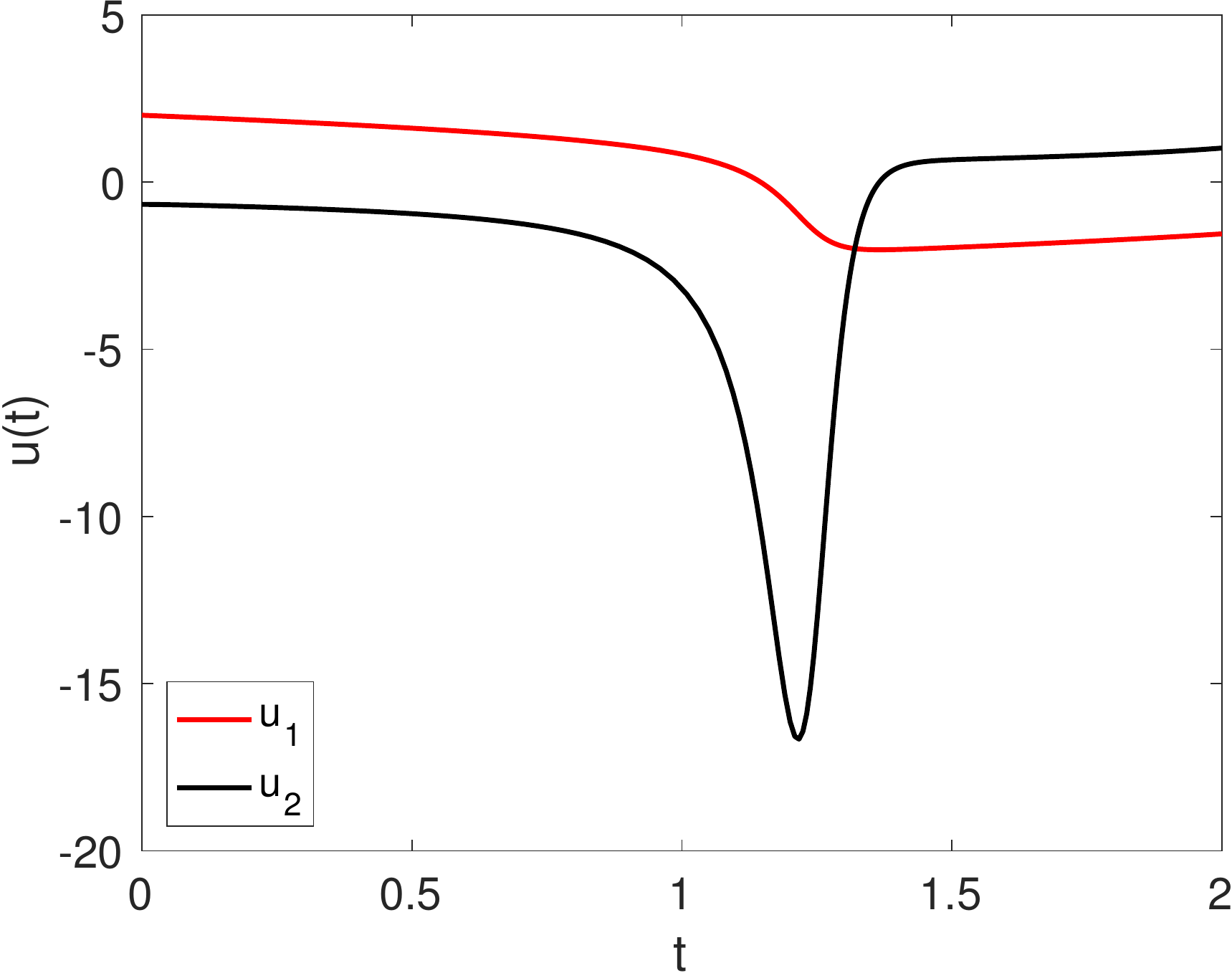}
        \caption{Van der Pol}
        \label{fig:vdpSolutionExample}
    \end{subfigure}
    \begin{subfigure}[b]{0.475\textwidth}
        \includegraphics[width=\textwidth]{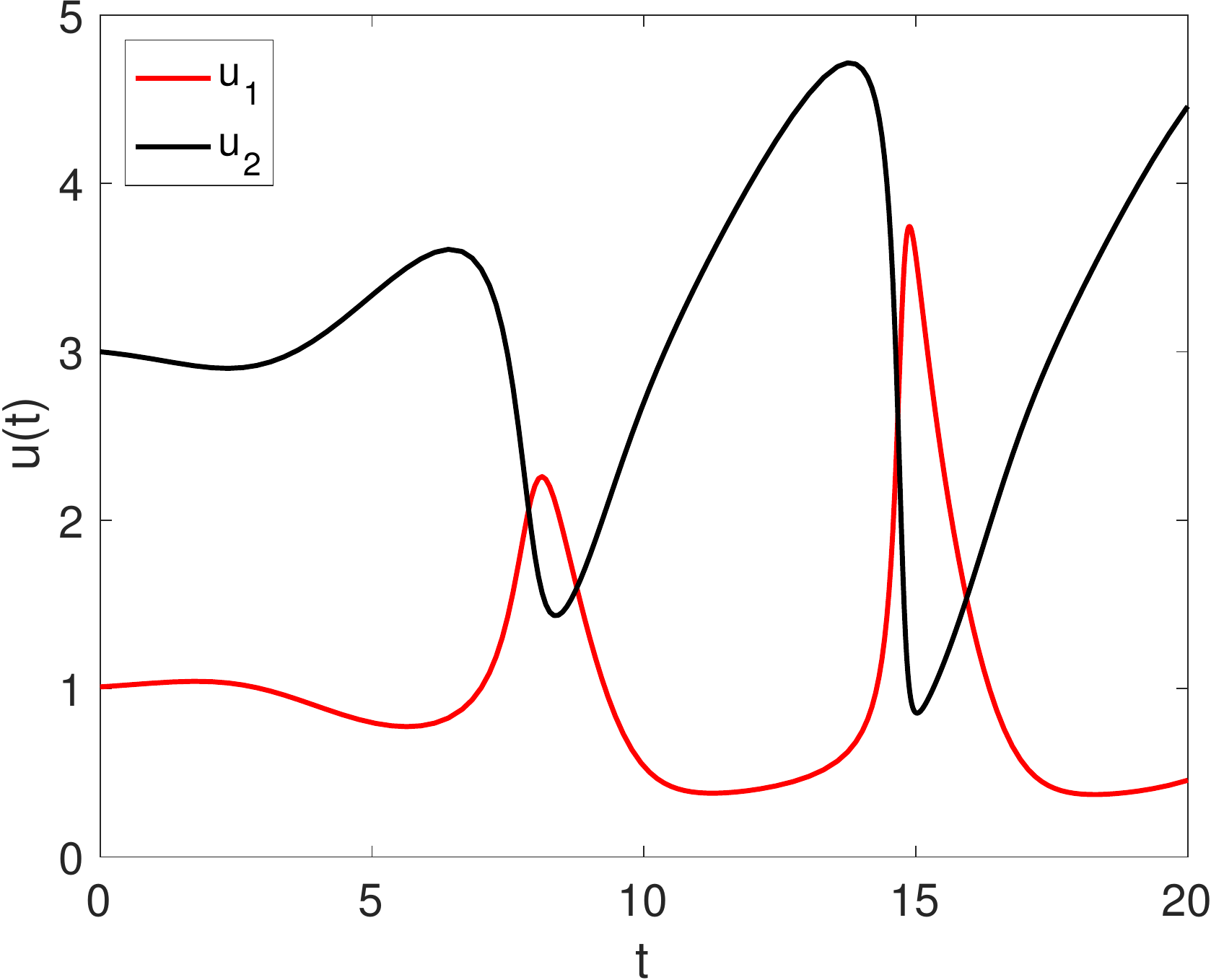}
        \caption{Brusselator}
        \label{fig:brusselatorSolutionExample}
    \end{subfigure} 
    \caption[Van der Pol and Brusselator Solution]{Van der Pol and Brusselator Solution.} 
    \label{fig:odeExampleSolution}
\end{figure}

We first begin with a relatively simple ODEs taken from \cite{iserles2008}: Van der Pol equation on  $0 \leq t \leq 2$  with $\epsilon = 0.1$ (Figure~\ref{fig:vdpSolutionExample})

\begin{equation}
\begin{aligned}
    u_1' &= u_2,  &  u_1(0) &= 2, \\
    u_2' &= \frac{1}{\epsilon} (1- u_1^2)u_2 - u_1, &  u_2(0) &= -0.6654321.
\end{aligned}
\label{eq:vdpProblem}
\end{equation}

\noindent Next we have another relatively simple ODEs: the Brusselator~\cite{HairerNorsettWanner93} on  $ 0 \leq t \leq 20$ (Figure~\ref{fig:brusselatorSolutionExample})

\begin{equation}
\begin{aligned}
    u_1' &= 1 + u_1^2u_2 - 4u_1 & u_1(0) &= 1.01, \\
    u_2' &= 3u_1 - u_1^2 u_2, & u_2(0) &= 3.
\end{aligned}
\label{eq:brusselatorProblem}
\end{equation}

\begin{figure}[ht!]
    \centering
    \begin{subfigure}[b]{0.475\textwidth}
        \includegraphics[width=\textwidth]{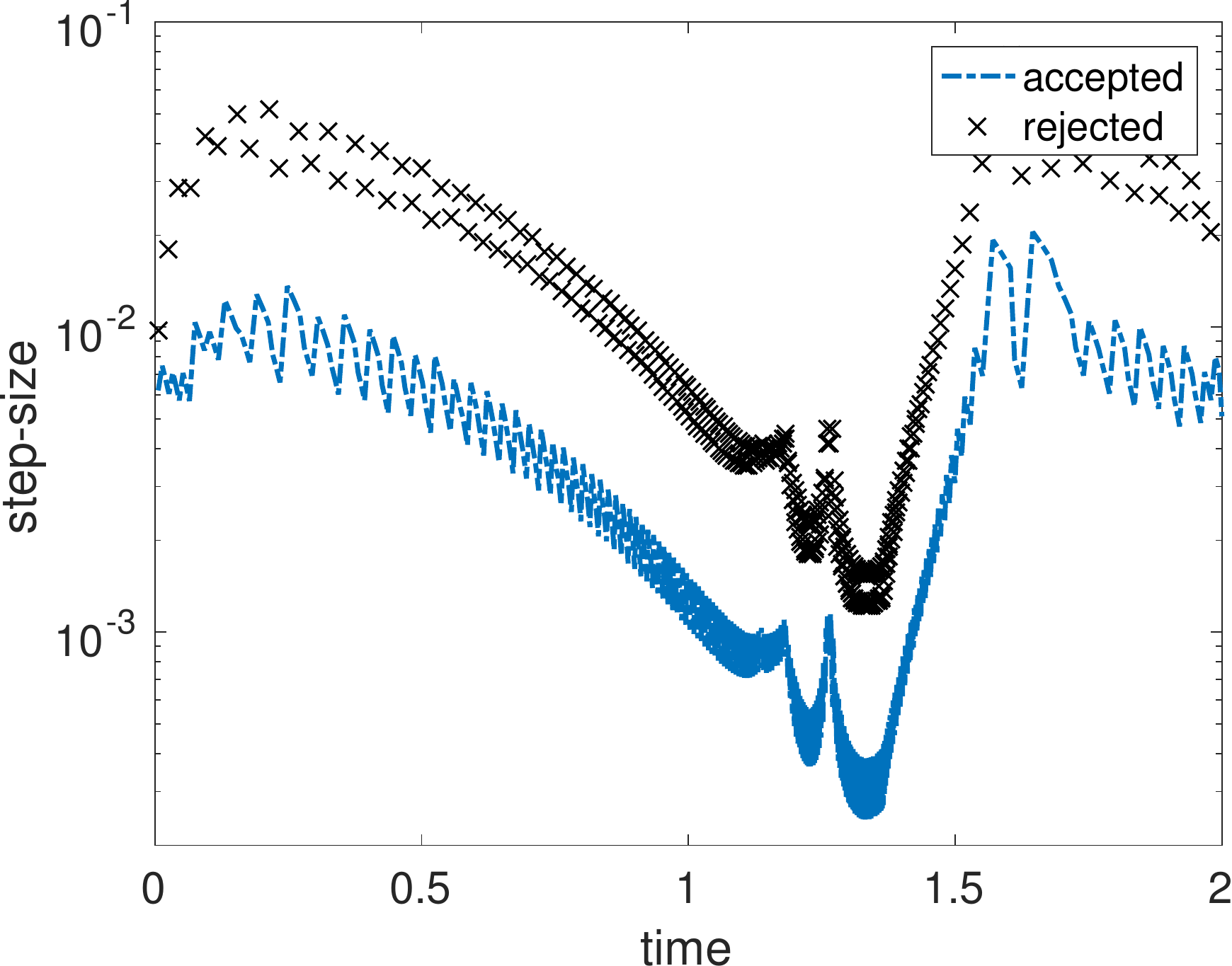}
        \caption{{I controller}}
        \label{fig:vdpStepSizeSecondOrderIcontroller}
    \end{subfigure} ~
    \begin{subfigure}[b]{0.475\textwidth}
       \includegraphics[width=\textwidth]{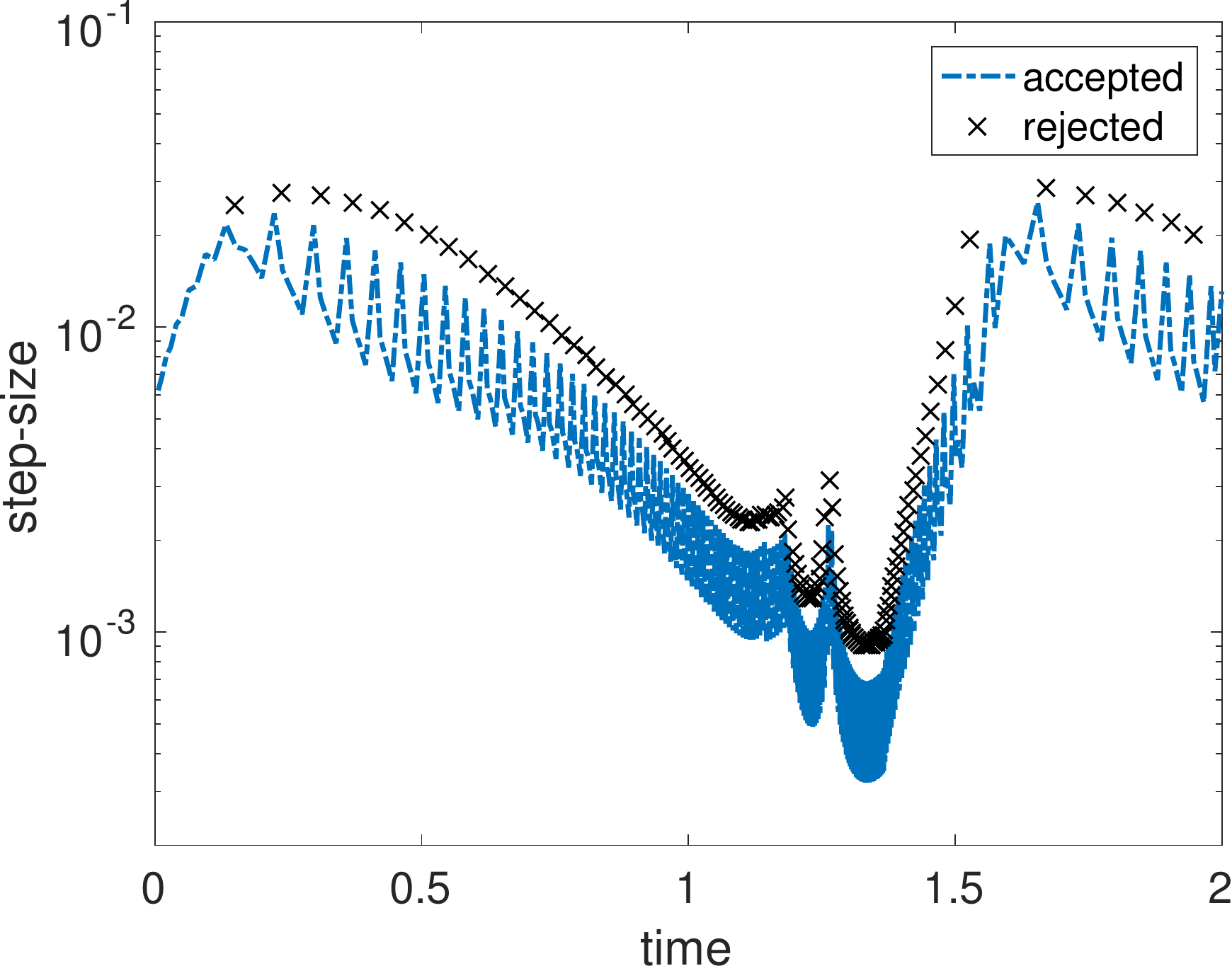}
        \caption{{PI controller}}
        \label{fig:vdpStepSizeSecondOrderPIController}
    \end{subfigure} \\
        \begin{subfigure}[b]{0.475\textwidth}
            \includegraphics[width=\textwidth]{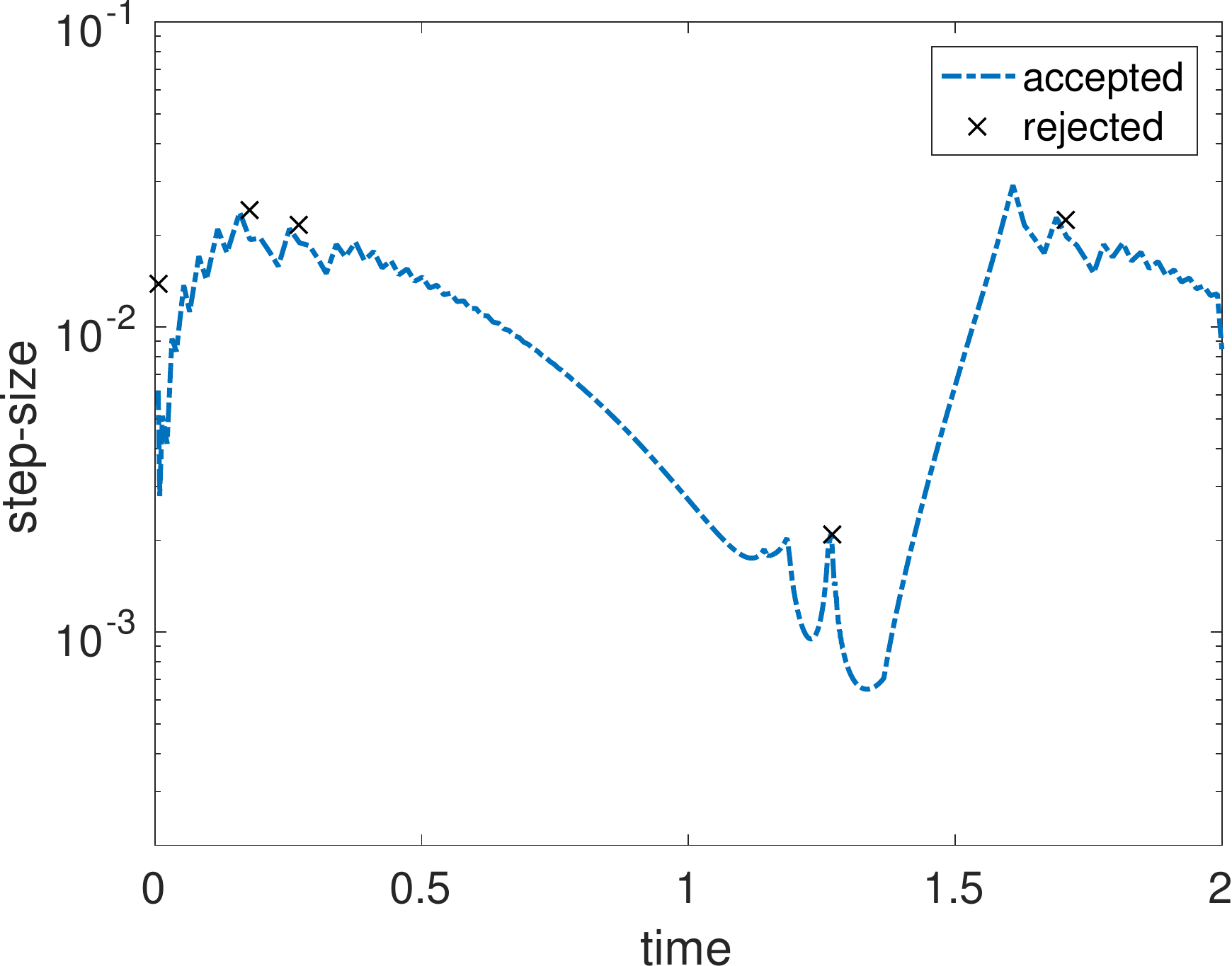}
            \caption{{PID controller}}
            \label{fig:vdpStepSizeSecondOrderPIDController}
        \end{subfigure}~
            \begin{subfigure}[b]{0.475\textwidth}
                \includegraphics[width=\textwidth]{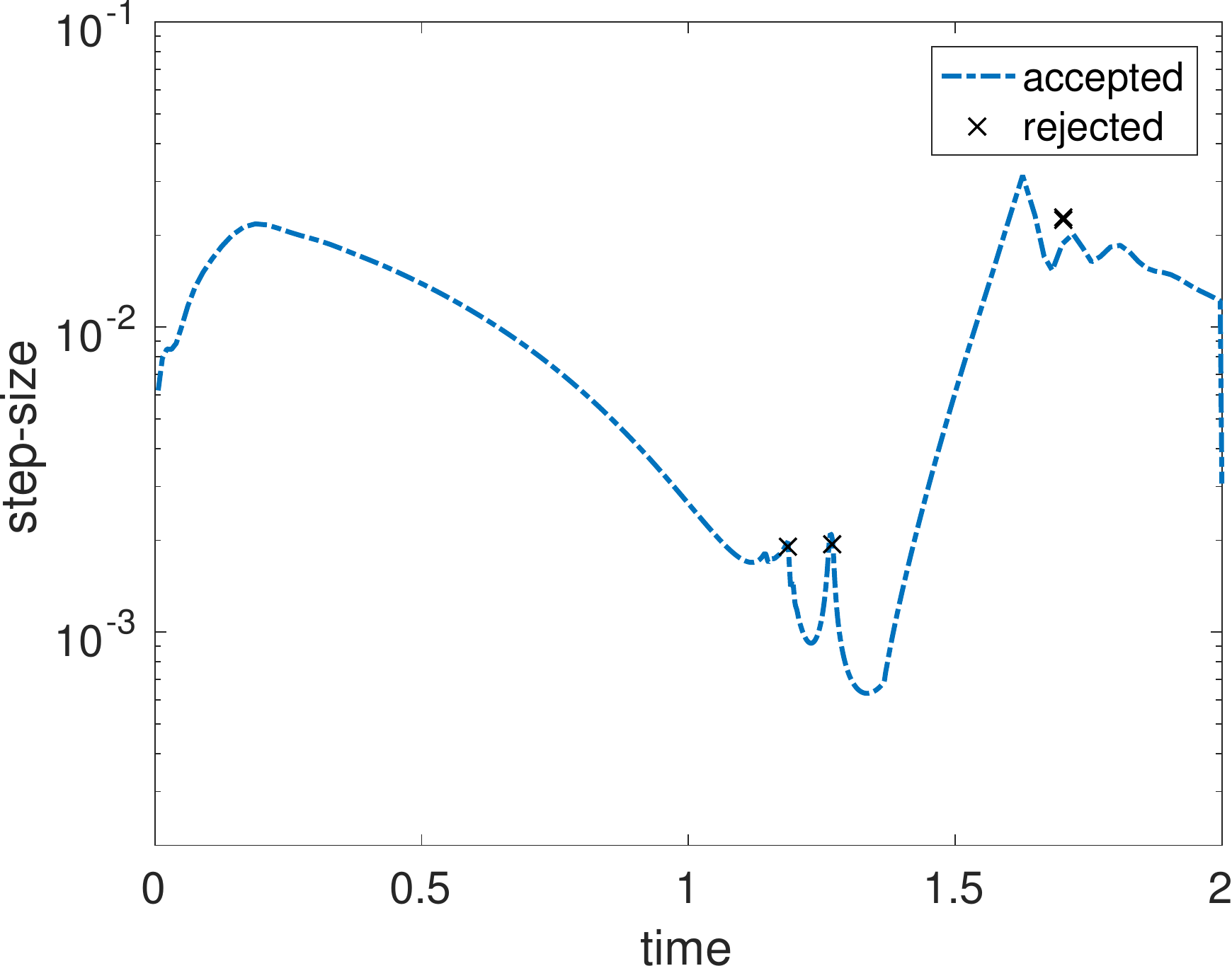}
                \caption{{Explicit Gustafsson controller}}
                \label{fig:vdpStepSizeSecondOrderExpGustController}
            \end{subfigure}
            \caption[Step-size comparison of SSPERK$(2,2)-\bt_2$ for VDP using the different adaptivity control algorithm at $\Rtol = \Atol = 1e{-4}$]{
            Step-size comparison of SSPERK$(2,2)-\bt_2$ for Van der Pol using the different adaptivity control algorithm at $\Rtol = \Atol = 1e{-4}$. Using the {I controller} (Figure~\ref{fig:vdpStepSizeSecondOrderIcontroller}): 1982 steps (495 rejected steps) with $L_2$-error = 4.06e-5. Using the {PI controller} (Figure~\ref{fig:vdpStepSizeSecondOrderPIController}): 1270 steps (210 rejected steps) with $L_2$-error = 1.09e-4. Using the {PID controller} (Figure~\ref{fig:vdpStepSizeSecondOrderPIDController}): 753 steps (17 rejected steps) with $L_2$-error = 1.59e-4. Using the {Explicit Gustafsson controller} (Figure~\ref{fig:vdpStepSizeSecondOrderExpGustController}): 795 steps (38 rejected steps) with $L_2$-error = 1.53e-4.} 
    \label{fig:vdpStepSizeSecondOrder}
\end{figure}

\begin{figure}[ht!]
    \centering
    \begin{subfigure}[b]{0.475\textwidth}
        \includegraphics[width=\textwidth]{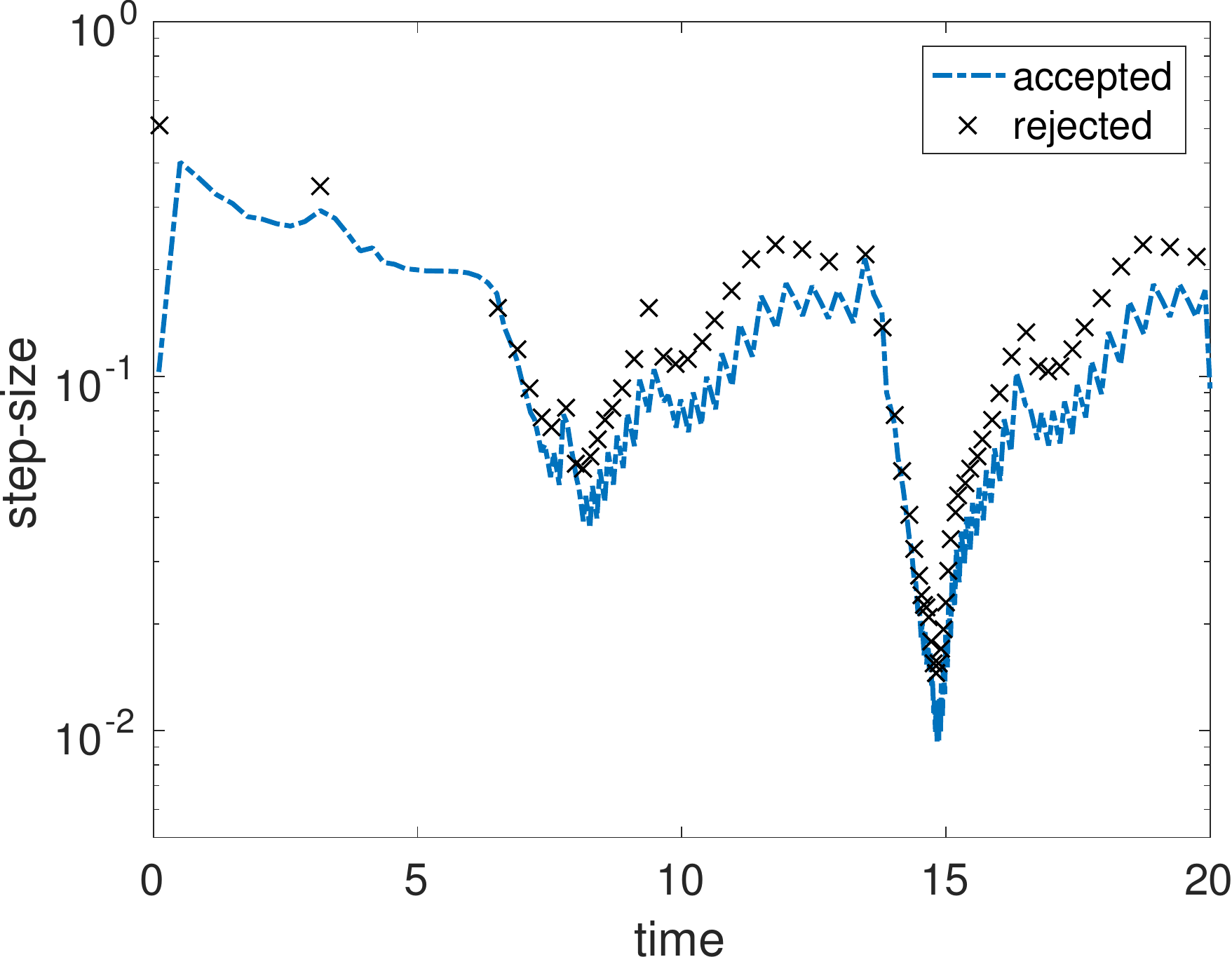}
        \caption{{I controller}}
        \label{fig:brusselatorStepSizeThirdOrderIcontroller}
    \end{subfigure} ~
    \begin{subfigure}[b]{0.475\textwidth}
       \includegraphics[width=\textwidth]{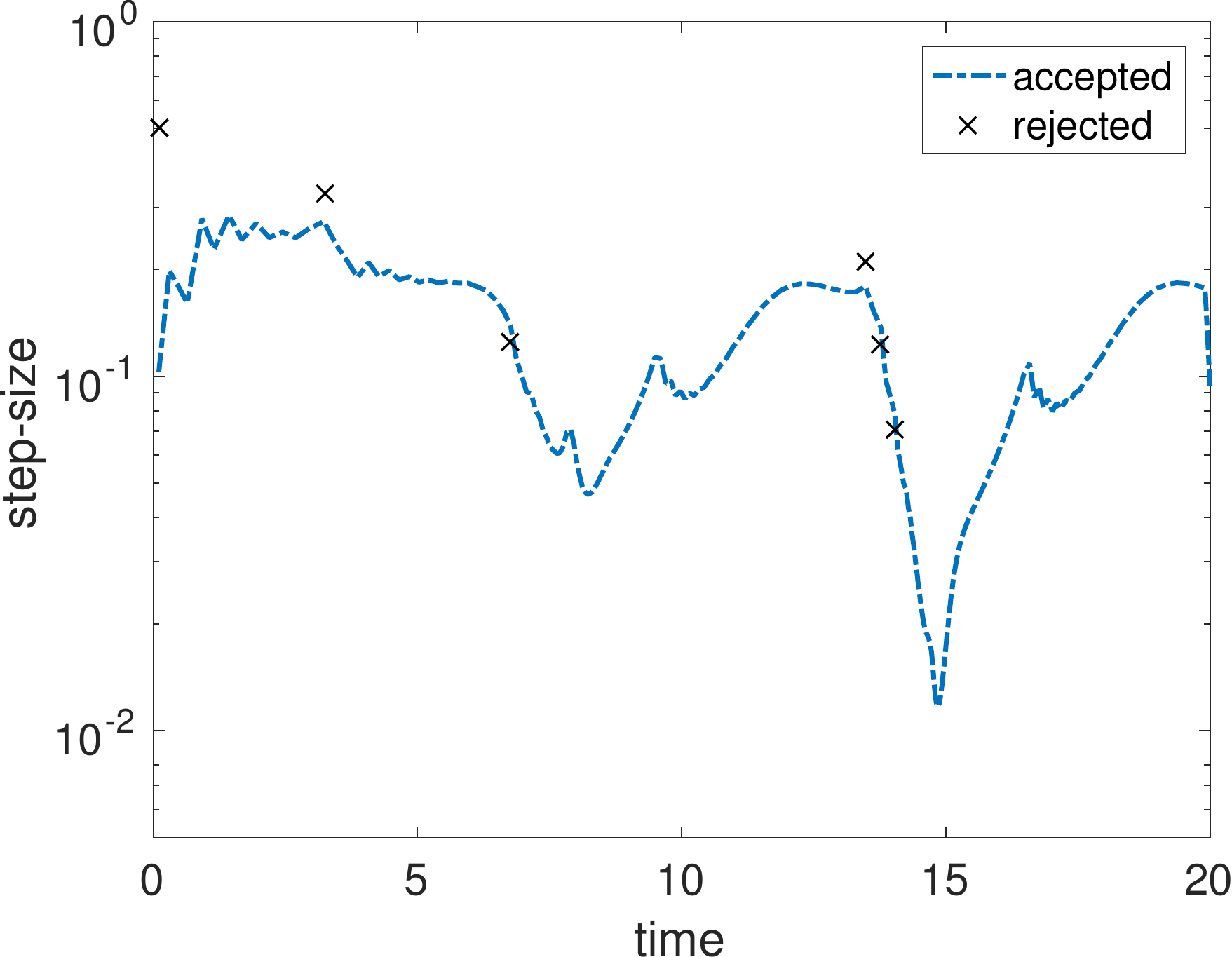}
        \caption{{PI controller}}
        \label{fig:brusselatorStepSizeThirdOrderPIController}
    \end{subfigure} \\
        \begin{subfigure}[b]{0.475\textwidth}
            \includegraphics[width=\textwidth]{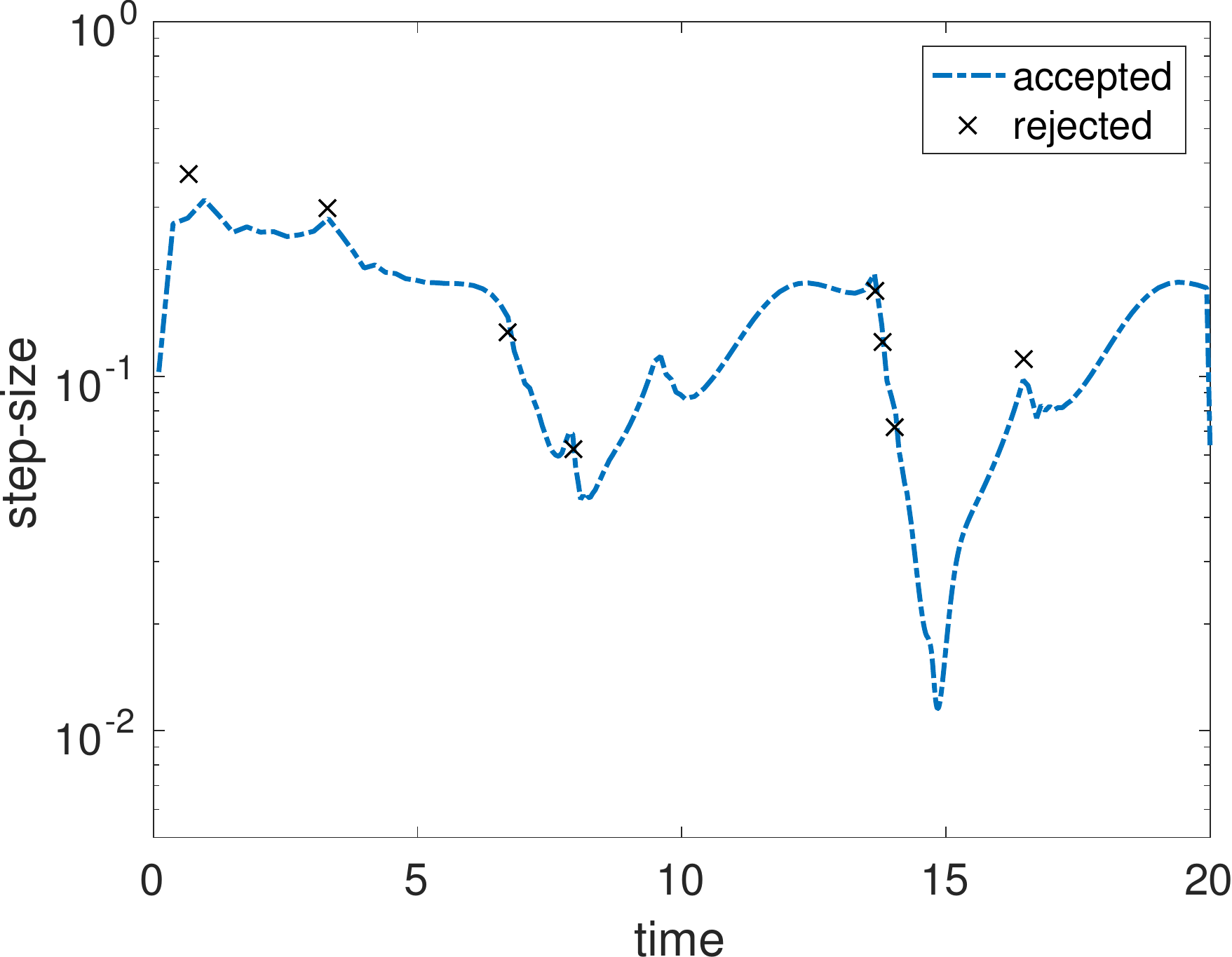}
            \caption{{PID controller}}
            \label{fig:brusselatorStepSizeThirdOrderPIDController}
        \end{subfigure}~
            \begin{subfigure}[b]{0.475\textwidth}
                \includegraphics[width=\textwidth]{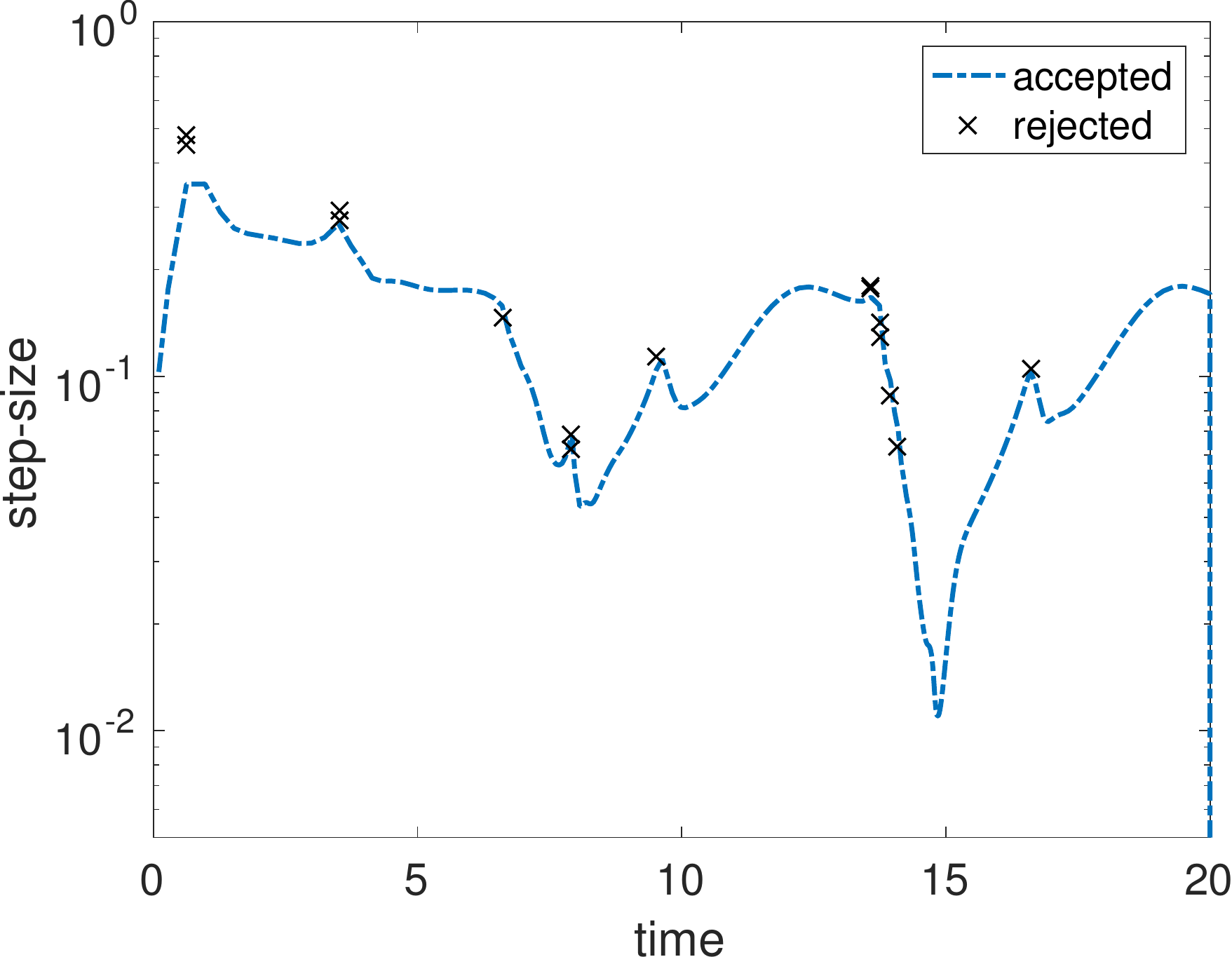}
                \caption{{Explicit Gustafsson controller}}
                \label{fig:brusselatorStepSizeThirdOrderExpGustController}
            \end{subfigure}
            \caption[Step-size comparison of SSPERK(3,3) for Brusselator using the different adaptivity control algorithm at $\Rtol = \Atol = 1e{-4}$]{Step-size for SSPERK(3,3) integrating Brusselator using the different adaptivity control algorithm at $\Rtol = \Atol = 1e{-4}$. Using the {I controller} (Figure~\ref{fig:brusselatorStepSizeThirdOrderIcontroller}): 419 steps (103 rejected steps) with $L_2$-error =2.7670e{-05}. Using the {PI controller} (Figure~\ref{fig:brusselatorStepSizeThirdOrderPIController}): 312 steps (17 rejected steps) with $L_2$-error = 3.2833e{-05}. Using the {PID controller} (Figure~\ref{fig:brusselatorStepSizeThirdOrderPIDController}): 305 steps (17 rejected steps) with $L_2$-error = 3.1775e{-05}. Using the {Explicit Gustafsson controller} (Figure~\ref{fig:brusselatorStepSizeThirdOrderExpGustController}): 332 steps (35 rejected steps) with $L_2$-error = 3.1086e-05.} 
    \label{fig:brusselatorStepSizeThirdOrder}
\end{figure}

\noindent We can see in Figure~\ref{fig:vdpStepSizeSecondOrder} the performance of the four control algorithms using a new SSPERK$(2,2)-\bt_2$ embedded pair on the Van der Pol test problem with $\Rtol = \Atol = tol = 1e{-4}$ prescribed. This figure shows all the step-sizes used; the accepted as well as the rejected ones (properly marked). We can see the {I Controller} often requires at least two attempts before proposing a successful optimal step-size.  It is more evident from Figure~\ref{fig:vdpStepSizeSecondOrderIcontroller} that the rejected step-size is far larger than the optimal accepted step-size throughout the integration. Furthermore, the average accepted step-size for this controller is $1.908e{-3}$.

The {PI Controller} takes fewer steps and requires less work than the {I Controller}; it has an average accepted step-size of $2.71e{-3}$.
In Figure~\ref{fig:vdpStepSizeSecondOrderPIController} we see less rejected steps. Moreover, the rejected step-sizes are closer to the accepted step-sizes than the I Controller. However this controller is inferior to the {PID Controller} which has an average accepted step-size of $3.937e{-3}$. The two previous controllers suggested optimal step-sizes are more sporadic than the {PID Controller} and {Explicit Gustafsson controller}. The average accepted step-size for Explicit Gustafsson controller is $3.839e{-3}$. In Figure~\ref{fig:vdpStepSizeSecondOrderPIDController} and Figure~\ref{fig:vdpStepSizeSecondOrderExpGustController}, we see the two results are much smoother. Furthermore, Figure~\ref{fig:vdpStepSizeSecondOrderPIDController} and Figure~\ref{fig:vdpStepSizeSecondOrderExpGustController} show far less total steps (and rejected steps), meaning fewer function evaluations needed to achieve the prescribed tolerance. In the vast majority of the numerical tests, the {PID Controller} performed significantly better. 


Through the numerous experiments, the PID and Explicit Gustafsson controller performed better than the I and PI controller. On some occasions, we noticed that the Explicit Gustafsson controller needed a second attempt to successfully predict an acceptable optimal step-size; the PID controller however, did not suffer from this limitation~(This is more evident in the third-order result for the long time Brusselator problem in Figure~\ref{fig:brusselatorStepSizeThirdOrder}). 

For the remaining comparisons, we show only the results using the {PID controller}. 
If we plot on a log-log scale, the ``work'' done by the method, which equates to the number of function evaluations calls and is measured as the number of stages~($s$) of the method multiplied by the total number of time-steps need to successfully integrate to the final time, versus the maximal global error~(``precision'') we obtain the a work-precision diagram. In what follows next, we describe the test problems and provide the associated work-precision results
; we use WENO5~\cite{SO88,jiang1996}, a fifth-order accurate method for the spatial discretization of the PDEs. 

\subsection{Numerial Performance}
\subsubsection{Numerical Performance: PDE Tests}
We consider two examples of the one-dimensional hyperbolic conservation law
\begin{equation}
 u_t + f(u)_x = 0
 \label{eq:scalarPDE}
\end{equation}

\noindent using WENO5~\cite{SO88,jiang1996} for spatial discretization. The first is a scalar linear advection equation $f(u) = u$ subject to periodic boundary conditions, with a square wave initial condition on $ x \in [-1, 1]$ and the integration interval $ 0 \leq t \leq 0.20$.

We then consider the more complex one-dimensional Euler Equations~\cite{laney1998computational,toro2009riemann}

\begin{equation} \label{eq:euler}
   \frac{\partial \bf{q}}{\partial t} + \nabla \cdot \bf{F} = 0 
\end{equation}

\noindent where the conserved variables, $\bf{q}$, and the flux, $\bf{F}$, are given as

\begin{equation*}
\bf{q} = \begin{bmatrix} \rho \\ \rho u \\ E \end{bmatrix} \quad
\bf{F} = \begin{bmatrix} \rho u \\ \rho u^2 + p \\ (E + p)u \end{bmatrix}
\end{equation*}

\noindent and the equations are closed by the ideal gas law as $$ p = (\gamma - 1) \left( E - \frac{1}{2} \rho u^2 \right), \quad c_s = \sqrt{\frac{\gamma p}{\rho}}$$
where $c_s$ represents the speed of sound and $\gamma$ is a fluid dependent constant. We take $\gamma =7/5$ for air in typical atmospheric conditions. We consider the Sod Shock tube~\cite{leveque2002finite,leveque1992numerical,Hesthaven2017} with initial conditions

\begin{align}
\rho(x,0) = 
  \begin{cases}
   1.0       & \quad x < 0.5 \\
   0.125       & \quad x \geq 0.5 \\
  \end{cases}
\qquad 
\rho u(x,0) = 0, 
\qquad
E(x,0) = \frac{1}{\gamma - 1}
  \begin{cases}
   1.0       & \quad x < 0.5 \\
   0.1       & \quad x \geq 0.5 \\
  \end{cases}
\end{align}

These also serve as boundary conditions since any disturbance is assumed not to reach the boundaries of the computational domain taken as $x \in [0, 1]$. The problem is integrated to a final time $t = 0.2$.


\subsubsection{Performance Results}

In what follows, we show the comparative performance on each test problems through the precision-work diagram~\cite{HairerNorsettWanner93}. Each result is compared against a reference method for easy comparison. The reference methods are: for second order results, SSPERK$(2,2)-\bt_2$; for third-order results, BogackiShampine32; and finally for fourth-order results, DormanndPrince54.
\begin{figure}[ht!]
    \centering
    \begin{subfigure}[b]{0.475\textwidth}
        \includegraphics[width=\textwidth]{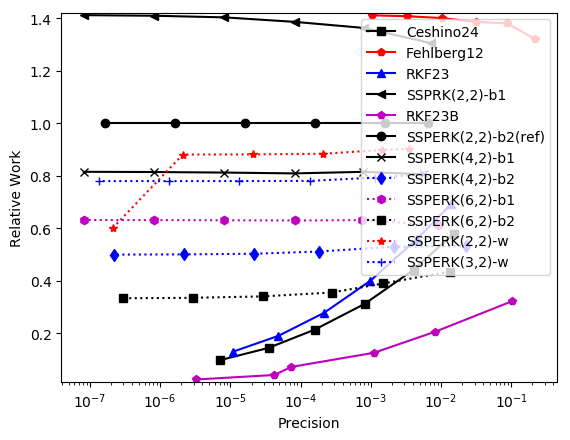}
        \caption{VDP}
        \label{fig:secondOrderVDP}
    \end{subfigure} ~
    \begin{subfigure}[b]{0.475\textwidth}
        \includegraphics[width=\textwidth]{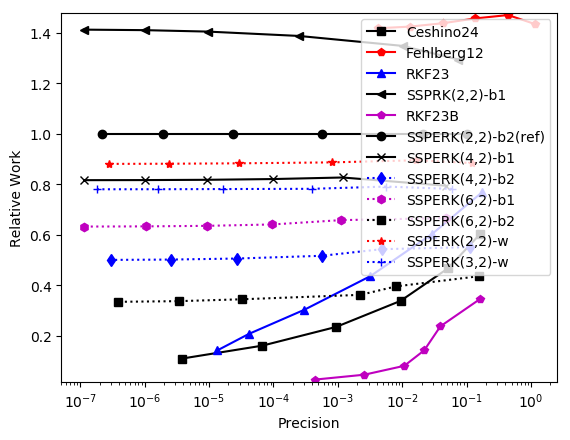}
        \caption{Brusselator}
        \label{fig:secondOrderBrusselator}
    \end{subfigure} \\
        \begin{subfigure}[b]{0.475\textwidth}
        \includegraphics[width=\textwidth]{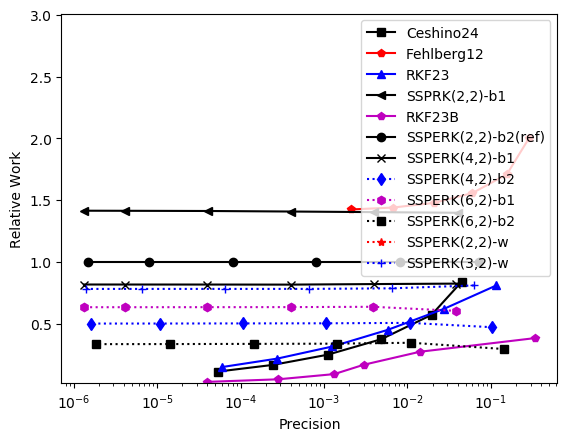}
            \caption{Advection}
            \label{fig:secondOrderAdvection}
        \end{subfigure}~
            \begin{subfigure}[b]{0.475\textwidth}
                \includegraphics[width=\textwidth]{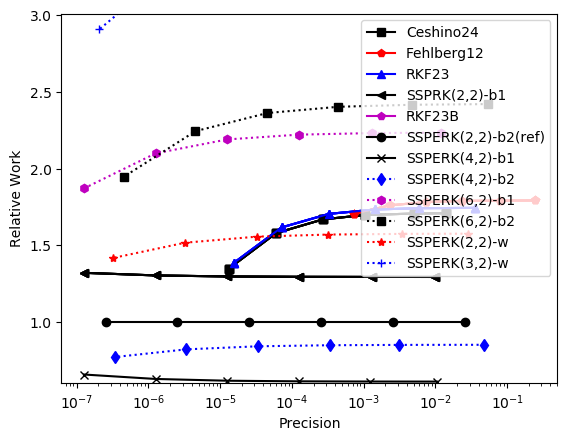}
                \caption{Euler}
                \label{fig:secondOrderEuler}
            \end{subfigure}
    \caption[Second order Precision-work diagram using {PID controller}]{Second order Precision-work diagram using PID controller. The ordinate is the relative work (number of function evaluations) compared against the SSPERK22-b2 pair. The abscissa is the global error at the endpoint of integration (the precision).} 
    \label{fig:PrecisionWorkSecondOrder}
\end{figure}

In Figure~\ref{fig:PrecisionWorkSecondOrder} we see the new SSPERK pairs are able to obtain a global error that is very close to the prescribed tolerance whereas the previous second order methods~(RKF23, Ceshino24, Felhber12) fall short. {Although it appears that these second-order methods~(RKF23, Ceshino24, Felhber12) are less costly than the newly developed SSPERK pairs, these methods fail to obtain a global error better than $1e{-5}$ even with a restrictive tolerance of $1e{-7}$ for any of the four problems.}
 
\begin{figure}[ht!]
    \centering
    \begin{subfigure}[b]{0.475\textwidth}
        \includegraphics[width=\textwidth]{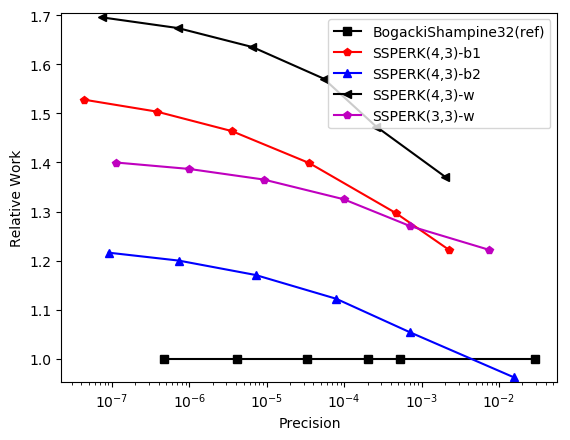}
        \caption{VDP}
        \label{fig:thirdOrderVDP}
    \end{subfigure} ~
    \begin{subfigure}[b]{0.475\textwidth}
        \includegraphics[width=\textwidth]{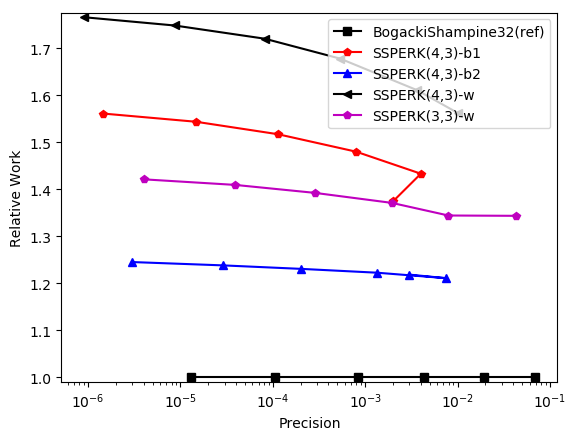}
        \caption{Brusselator}
        \label{fig:thirdOrderBrusselator}
    \end{subfigure} \\
        \begin{subfigure}[b]{0.475\textwidth}
            \includegraphics[width=\textwidth]{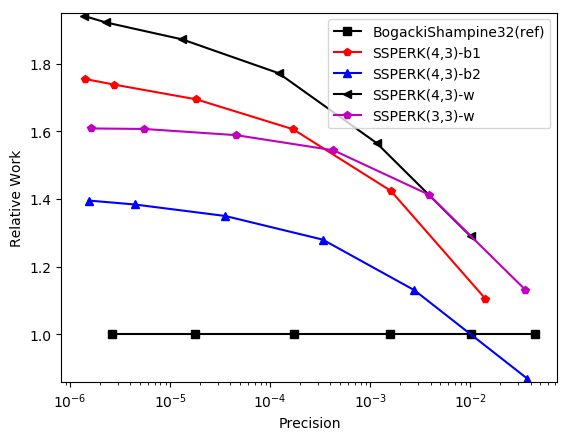}
            \caption{Advection}
            \label{fig:thirdOrderAdvection}
        \end{subfigure}~
            \begin{subfigure}[b]{0.475\textwidth}
                \includegraphics[width=\textwidth]{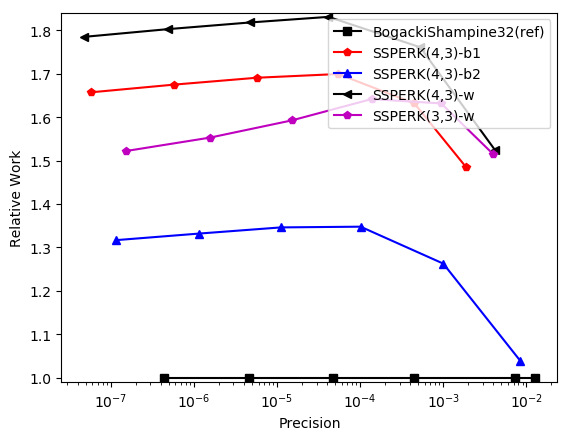}
                \caption{Euler}
                \label{fig:thirdOrderEuler}
            \end{subfigure}
    \caption[Third order Precision-work diagram using {PID controller}]{Third order Precision-work diagram using PID controller. The ordinate is the relative work (number of function evaluations) compared against the BogackiShampine32 pair. The abscissa is the global error at the endpoint of integration (the precision).} 
    \label{fig:PrecisionWorkThirdOrder}
\end{figure}

Figure~\ref{fig:PrecisionWorkThirdOrder} shows the relative work-precision results for the third-order methods for each of the four problems. 
To no surprise, the numerical experiments proved that BogackiShampine(3,2) has better error measurements than the optimal SSPERK pairs. {The results show the new SSPERK pairs are more expensive when compared to the methods by Bogacki and Shampine (BogackiShampine32)}. 
It is worth noting this 4 stage method has all positive coefficient RK(3,2) pair but is not SSP. For the hyperbolic problems, the newly constructed method SSPERK$(4,3)-\bt_2$ performs relatively similar to the BogackiShampine(3,2) pairs for $tol \geq 1e-2$ and only 30\%-40\% more costly otherwise. At this tolerance, it is advantageous to use SSPERK$(4,3)-\bt_2$ over BogackiShampine32 since the performance is relatively the same and the SSPERK$(4,3)-\bt_2$ pair offers the benefits of strong stability preservation. 

\begin{figure}[ht!]
    \centering
    \begin{subfigure}[b]{0.475\textwidth}
        \includegraphics[width=\textwidth]{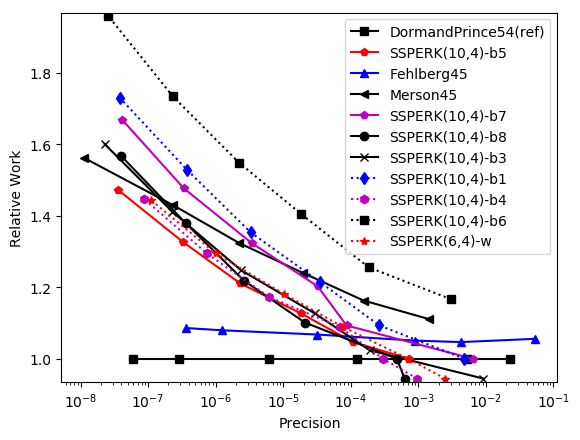}
        \caption{VDP}
        \label{fig:fourthOrderVDP}
    \end{subfigure} ~
    \begin{subfigure}[b]{0.475\textwidth}
        \includegraphics[width=\textwidth]{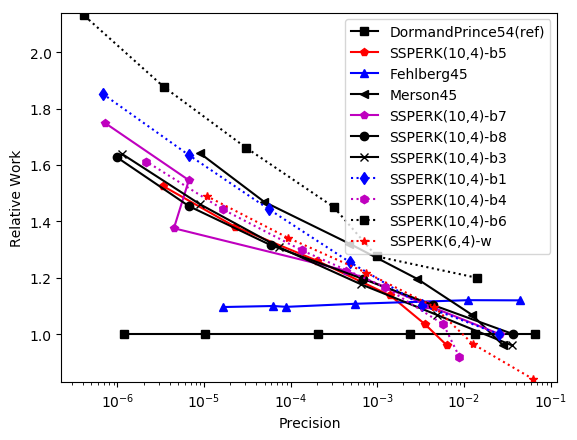}
        \caption{Brusselator}
        \label{fig:fourthOrderBrusselator}
    \end{subfigure} \\
        \begin{subfigure}[b]{0.475\textwidth}
            \includegraphics[width=\textwidth]{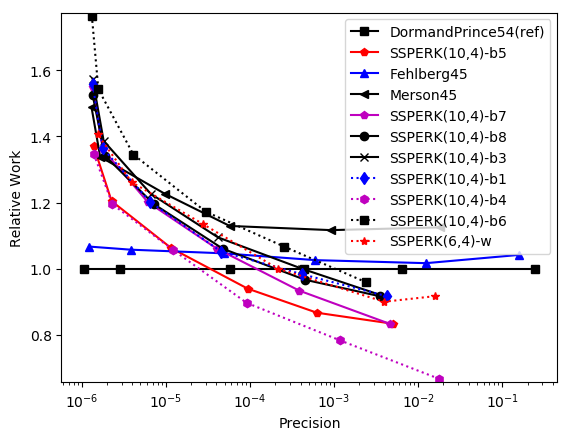}
            \caption{Advection}
            \label{fig:fourthOrderAdvection}
        \end{subfigure}~
            \begin{subfigure}[b]{0.475\textwidth}
                \includegraphics[width=\textwidth]{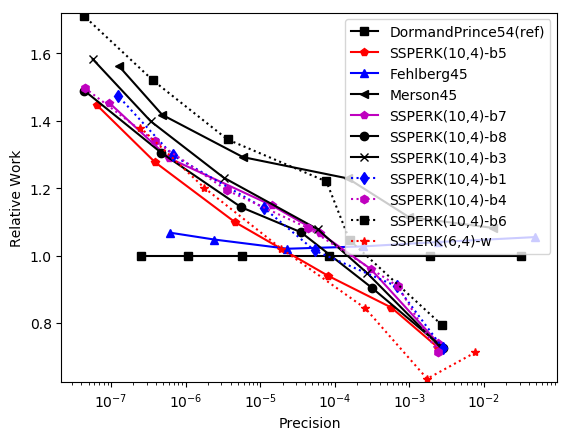}
                \caption{Euler}
                \label{fig:fourthOrderEuler}
            \end{subfigure}
    \caption[Fourth order Precision-work diagram using {PID controller}]{Fourth order Precision-Work diagram using PID controller. The ordinate is the relative work (number of function evaluations) compared against the DormandPrince54 pair. The abscissa is the global error at the endpoint of integration (the precision).} 
    \label{fig:PrecisionWorkFourthOrder}
\end{figure}

Similarly, Figure~\ref{fig:PrecisionWorkFourthOrder} shows the relative work-precision results for the fourth-order methods for each of the test problems. The best performing pairs for the SSPERK$(10,4)$ are shown as well as the SSPERK$(6,4)-\wt$. 
These new embedded pairs are compared to the five-stage Merson45, six-stage Fehlberg45, five-stage Zonneveld43 and the seven-stage DormandPrince54. The DormanPrince54 method is used in MATLAB's \emph{ODE45} implementation. With at least one negative coefficient, none of the fourth-order literature embedded pairs are SSP. 
The new SSPERK methods are more efficient than the other methods for large tolerance $tol \geq 1e-4$ for hyperbolic problems~(Figures \ref{fig:fourthOrderAdvection}~and~\ref{fig:fourthOrderEuler}). 
At tolerance $tol \geq 1e-4$, larger step-size can be taken and may violate the stability requirement. However the newly constructed SSPERK pairs perform best in these regions. At a very restrictive tolerance ($tol \leq 1e-5$) the optimal step-sizes are small enough that the methods don't encounter any stability issues and SSPERK pairs are no longer favorable. 

We note that the SSPERK$(10,4)-\bt_2$ performed the worst for all problems. In all the numerical studies, this method overestimates the local error and returns a really small step-size causing the method to do an incredibly large amount of work to integrate the problem to the final time. At the final time, the global error (the precision) is much lower than the prescribed tolerance (i.e. $1e{-12}$ vs. $1e{-7}$). A practical error estimate ``ensures that the step-sizes are sufficiently large to avoid unnecessary computational work"~\cite{HairerNorsettWanner93}. SSPERK$(10,4)-\bt_2$ is not practical since it fails to do this. 

\section{Remarks}\label{sec:discussion}

SSP methods are used heavily for numerically integrating IVP resulting from discretized hyperbolic problems since they offer properties not present in non-SSP methods. Modern robust IVP solvers include many important features such as error estimation and automatic step-size control~\cite{C00,Carpenter2000,ShampineReichelt97}. Most of these important features have not yet been developed for existing higher-order optimal explicit SSP methods. The current work provides these important features for the optimal SSP methods. The numerical results provide evidence of the robustness of these pairs. 

Although the existing pairs are not SSP, we note that using a perturbation technique from~\cite{ketchesonHiguerasDWRK} can provide a positive SSP coefficient and these features can be used in the usual way at the expense of an additive routine. 

 {The total work, i.e.,computational cost, scales with the stage number of the method and the total number of steps required to integrate the problem to $\tfinal$: large stage number requires more function evaluations per single step computation but often requires fewer total number of steps. A method with large number of stages is often equipped with larger stability regions. For SSP methods,  higher stage number leads to higher SSP coefficient $\sspcoef$. In the study we see that pairs with large stage counts are very robust and perform less work per degree of accuracy: increasing the number of stages can reduce the total work. 
 These factors may indicate that having a large number of stages, characterized by a large stability region and also large SSP coefficients, is critical for simulations where the total time to obtain a sufficiently accurate result is paramount.
 
We can conclude, depending on the problem being integrated and length of integration, the new SSPERK embedded pairs are very effective, practical, and robust. Considering the time to accurate and stable solution, we make the following recommendations. For second-order methods, we recommend the pairs SSPERK$(4,2)-\bt_2$,  SSPRK$(6,2)-\bt_2$ and due to its low stage count, the SSPERK$(3,2)-\wt$. For third-order methods, SSPERK$(4,3)-\bt_2$ performed closest to the method of Bogacki and Shampine. Due to its popularity, the pair for SSPERK$(3,3)$ is also recommended. {Finally apart from the  inpractical pair SSPERK$(10,4)-\bt_2$, the results from the remaining SSPERK$(10,4)$ are relatively similar. Though we recommend the pairs SSPERK$(10,4)-\bt_1$, SSPERK$(10,4)-\bt_3$, SSPERK$(10,4)-\bt_5$, and SSPERK$(10,4)-\bt_8$.}

\section*{Acknowledgement}
The work of Sidafa Conde and John Shadid was partially supported by the U.S. Department of Energy , Office of Science, Office of Applied Scientific Computing Research. Sandia National Laboratories is a multi-mission laboratory managed and operated by National Technology and Engineering Solutions of Sandia, LLC., a wholly owned subsidiary of Honeywell International, Inc., for the U.S. Department of Energys National Nuclear Security Administration under contract DE-NA0003525. The views expressed in the article do not necessarily represent the views of the U.S. Department of Energy or the United States Government. The work of Imre Fekete was supported by the Hungarian Scientific Research Fund OTKA under grants No. 112157 and the project has been supported by the European Union, co-financed by the European Social Fund (EFOP-3.6.3-VEKOP-16-2017-00002). In addition, this work used the Extreme Science and Engineering Discovery Environment (XSEDE), which is supported by National Science Foundation grant number ACI-1548562~\cite{xsederesources}. This work used the Extreme Science and Engineering Discovery Environment (XSEDE) Stampede at the Texas Advanced Computing Center (TACC) through allocation DMS170002.

The authors thank Sigal Gottlieb and David Ketcheson for the stimulating discussions on the topic and their helpful remarks.


\bibliographystyle{plain}
\bibliography{embedded_ssp_methods.bib}
\end{document}